# Classification of Collatz infinite sequences

Raouf Rajab
raouf.rajab@enig.rnu.tn

**Abstract**

In the present paper, we are interested in classifying of Collatz sequences on based to the different behavior of these sequences when their lengths tend to infinity. A Collatz infinite sequence can be defined as an infinite ordered set of positive integers such that the term of rank n is results of applying Collatz map n times to the first term. Such term can be expressed on the form $T^n(P)=A(P,n)P+B(P,n)$. When n tends to infinity, each function among the two partial coefficients denoted by $A(P,n)$ and $B(P,n)$ behaves in different ways. This allows us to determine all categories of Collatz infinite sequences. First, we carry out a classification of Collatz infinite sequences on based of the different possible limits of the two coefficients. In second time, we determine the different proportions of every class of the infinite sequences. Note that results obtained do not represent a proof of a Collatz conjecture but they have a strong relationship with this conjecture and it allows us to better understand the behavior of Collatz sequences when n tend to infinity.

**Résumé :** Dans le présent travail, on s'intéresse à la catégorisation des suites de Collatz. Une telle classification est basée sur le comportement des suites de Collatz lorsque ses longueurs tendent vers l'infini. On définit une suite infinie de Collatz comme étant un ensemble infini ordonné des entiers naturels non nuls tel que le terme d'ordre n de cette suite est obtenu par application de la fonction de Collatz n fois successive à l'entier de départ. Un tel terme s'écrits sous la forme $T^n(P)=A(P,n)P+B(P,n)$. Lorsqu'on fait tendre n vers l'infini les deux fonctions partielles $A(P,n)$ et $B(P,n)$ se comportent des différentes manières. En premier temps, on s'intéresse à la classification des suites de Collatz de longueurs infinies en se basant sur les différentes limites possibles des coefficients $A(P,n)$, $B(P,n)$ et de $T^n(P)$. Puis on détermine les différentes proportions relatives aux différentes suites infinies appartenant aux différentes catégories déjà définies. Noter que les résultats obtenus ne représentent pas une démonstration de la conjecture de Collatz mais ils sont en forte relation avec cette conjecture et nous permet de mieux comprendre les propriétés de convergences de ce genre des suites.

**Mots clés :** Suites infinies de Collatz ; Coefficients caractéristiques absolus; Suites infinies réalisables; Suites infinies non réalisables.

## 1. Introduction

La conjecture de Collatz largement connue sous le nom du problème 3n+1 ou aussi conjecture de Syracuse (on peut rencontrer plusieurs autres nominations) est l'un des problèmes ouverts les plus fameux en mathématiques. Elle fait l'objet des plusieurs recherches au cours de laquelle plusieurs pistes ont été explorées par un grand nombre des chercheurs dans l'objectif final est la démonstration de cette conjecture qu'on peut le résumer en deux problèmes [1]:

-Toutes les suites de Collatz atteignent 1 après un certain nombre fini des itérations successives.

-Le cycle (1,2) est le seul cycle stable par itération de la fonction de Collatz.

Cette conjecture est basée sur un processus itératif décrit par la fonction de Collatz ci-dessous qu'on la note T. Elle est définie de $\mathbb{N}^*$ dans $\mathbb{N}^*$ pour tout entier naturel non nul P comme suit [2]:

(1.1) $$T(P) = \frac{3^{\mathbf{i}_0(P)}}{2}P + \frac{1}{2}\mathbf{i}_0(P)$$

Avec $\mathbf{i}_0$ est l'indicateur de parité, il est définie comme une application de $\mathbb{N}$ dans {0,1} tel que pour tout entier naturel P:

(1.2) $$\mathbf{i}_0(P) = \begin{cases} 0 & \text{si } P \equiv 0 \pmod{2} \\ 1 & \text{si } P \equiv 1 \pmod{2} \end{cases}$$

Dans ce travail, on va adopter les notations suivantes

**Notation        1.1**
Pour tout entier naturel non nul P, on définit la suite de Collatz de longueur n comme suit [2]:

(1.3) $$Sy(P, n) = \left(P, T^1(P), T^2(P), T^2(P), \ldots, T^n(P)\right)$$

Si les suites considérées sont des longueurs infinies, on adopte la notation $Sy^\infty(P)$.

**Notation        1.2**
On fait associer à chaque suite de Collatz de premier terme P et de longueur n notée $Sy(P, n)$ un vecteur de parité qu'on le note $v(P, n)$ [3].

(1.4) $$Sy(P, n) \rightarrow v(P, n)$$

**Notation        1.3**
On désigne par ClZ l'ensemble constitué par toutes les suites de Collatz qui atteignent 1 après un nombre fini des itérations à partir d'un entier naturel nul de départ cet ensemble est définie comme suit:

(1.5) $$CLZ = \{Sy^\infty(P) \text{ tel que } T^n(P) = 1 \text{ et } (P, n) \in \mathbb{N}^* x \mathbb{N}^*\}$$

**Notation        1.4**
Le processus itératif de Collatz nous conduit à l'obtention d'une image d'ordre n de P qui s'écrit sous la forme :

(1.6) $$T^n(P) = A_n(P)P + B_n(P)$$

Ou $A_n(P)$ et $B_n(P)$ sont deux fonctions de $\mathbb{N}^*$ dans $\mathbb{Q}^*$ résultantes de l'application de la fonction de Collatz n fois successives sur P qu'on peut l'appeler les coefficients caractéristiques d'ordre n de la suite de Collatz de longueur n.

$A_n(P)$ Coefficient principal d'ordre n

$B_n(P)$ Coefficient secondaire ou d'ajustement d'ordre n.

**Remarque 1.1**

On peut adopter les notations $A(P, n)$ et $B(P, n)$ au lieu de $A_n(P)$ et $B_n(P)$ et vise versa.

On peut étudier les propriétés des convergences des suites de Collatz à partir de l'étude des comportements de deux coefficients (ou fonctions partielles) lorsque n tend vers l'infini.

Le nombre des itérations nécessaires pour qu'une suite donnée atteint le cycle (2,1) varie d'une suite à une autre mais les deux conditions qui sont remplies par toutes les suites de Collatz appartenant à l'ensemble **CLZ** sont les suivantes :

Les coefficients principaux des toutes les suites infinies de l'ensemble CLZ sont pratiquement nuls :

$$\lim_{n \to +\infty} A_n(P) = 0$$

Le deuxième coefficient tend vers le cycle (1,2) stable par la transformation de Collatz:

$$B_n(P) \to (1,2)$$

Pour plus de précision sur les comportements de ce coefficient, on peut distinguer deux différents cas :

Pour une suite qui atteint 2 pour la première fois après un nombre impair d'itérations:

$$\begin{cases} \lim_{k \to +\infty} B_{2k}(P) = 1 \\ \lim_{k \to +\infty} B_{2k-1}(P) = 2 \end{cases}$$

Pour une suite qui atteint 2 pour la première fois après un nombre pair d'itérations:

$$\begin{cases} \lim_{k \to +\infty} B_{2k-1}(P) = 1 \\ \lim_{k \to +\infty} B_{2k}(P) = 2 \end{cases}$$

Autrement, pour une suite de Collatz de premier terme P vérifiant la conjecture de Collatz nécessairement lorsque n tend vers l'infini le coefficient $A_n(P)$ tend vers 0 alors que le deuxième coefficient tend vers le cycle (1,2). Une telle suite va rejoindre donc la suite de Collatz de premier terme 2 après un certain nombre des itérations.

En général, on peut supposer qu'on a deux cas possibles pour la limite de $A_n(P)$ à l'infini et aussi deux cas possibles pour la limite de $B_n(P)$ donc si on fait combiner tous les cas possibles on obtient 4 possibilités pour les comportements des suites infinies.

Tableau1 : Les quatre cas possibles pour les comportements de $A_n(P)$ et de $B_n(P)$ à l'infini.

|  | $\lim_{n\to+\infty} A_n(P)$ est finie | $\lim_{n\to+\infty} A_n(P)$ est infinie |
|---|---|---|
| $\lim_{n\to+\infty} B_n(P)$ est finie | $(-,-)$ | $(+,-)$ |
| $\lim_{n\to+\infty} B_n(P)$ est infinie | $(-,+)$ | $(+,+)$ |

En réalité, ces quatre Cas ou situations composées ne sont pas suffisantes pour décrire entièrement les comportements des toutes les suites à l'infini, on peut distinguer d'autres cas qui sont liés à la convergence de $T^n(P)$, ou bien la comparaison entre le premier terme d'une suite et sa limite à l'infini , on peut faire la distinction entre les suites **infinies** qui admettent des coefficients principaux nuls et celles qui possèdent des coefficients principaux non nuls…exc.

Le présent travail ne porte pas sur la démonstration de cette conjecture mais on s'intéresse à la classification des suites de Collatz en se basant sur les différentes limites possibles de deux fonction partielles $A_n(P)$ et $B(P,n)$. Cette classification nous permet de faire la distinction entre les différentes catégories des suites infinies de Collatz et nous permet de mieux comprendre leurs comportements à l'infini. On fait combiner plusieurs propriétés liées aux comportements de deux fonctions $A_n(P)$ et $B(P,n)$ lorsque n tend vers l'infini, on peut remarquer que quelques combinaisons correspondent à des cas impossibles ou a des catégories vides. Puis on détermine la proportion des suites appartenant à une même catégorie dans l'ensemble de toutes les suites de Collatz de longueurs infinies. Une telle catégorisation est en forte relation avec la conjecture de Collatz et elle permet de faciliter l'étude de ce genre de suites.

## 2. Préliminaire

**Définition      2.1**

Soit P un entier naturel non nul, pour tout entier naturel n (n ≥ 0 ), on définit l'indicateur de parité d'ordre n de P qu'on le note $\mathbf{i}_n(P)$ comme suit :

(2.1) $$\mathbf{i}_n(P) = \begin{cases} 1 \text{ si } T^n(P) \text{ est impair} \\ 0 \text{ si non} \end{cases}$$

**Notation      2.1**

Le nombre des entiers impairs contenus dans une suite de Collatz $S^c(P,n)$ est noté $m_n(P)$ et il peut être exprimé en fonction des indicateurs de parité comme suit :

(2.2) $$m_n(P) = \sum_{k=1}^{k=n} \mathbf{i}_k(P)$$

### 2.1 Matrice complète d'ordre n de Collatz

On considère les deux ensembles suivants des suites de collatz de même longueur n :

$$\begin{cases} D1 = \big(Sy(1,n), Sy(3,n), Sy(5,n), \dots, Sy(2^{n+1}-1, n)\big) \\ D2 = \big(Sy(2,n), Sy(4,n), Sy(6,n), \dots, Sy(2^{n+1}, n)\big) \end{cases}$$

Chaque ensemble est représenté sous forme d'un tableau de ( n+1) colonnes et de $2^n$ lignes. La première colonne contient les $2^n$ premiers termes de différentes suites d'un même ensemble. Les deux tableaux construits sont appelés les matrcies complètes de Collatz d'ordre dimentionnel n.

|           | $T^1(P_i)$ |  |  | $T^n(P_i)$ |
|-----------|------------|--|--|------------|
| $2^{n+1}-1$ |          |  |  |            |
|           |            |  |  |            |
|           |            |  |  |            |
| $P_i$     |            |  |  |            |
|           |            |  |  |            |
| 5         |            |  |  |            |
| 3         |            |  |  |            |
| 1         |            |  |  |            |

$$\begin{cases} A_n(P) \\ B_n(P) \\ T^n(P) \end{cases} \longrightarrow \begin{cases} A_\infty(P) \\ B_\infty(P) \\ T^\infty(P) \end{cases}$$

Figure 1: Matrice complète de Collatz d'ordre dimentionnel n : $\mathbb{M}^T(1, n)$

On peut adopter la notation suivante pour une suite de collatz de longueur n contenue dans la matrcie considérée :

(2.3) $\qquad\qquad S^c(P, n) = \big(T^1(P), T^2(P), T^2(P), \dots, T^n(P)\big)$

Noter que P n'appartient pas à la suite $S^c(P, n)$.

On fait tendre n vers l'infini, on obtient les deux tableaux complets d'ordre dimensionnel infini :

Le premier tableau contient toutes les suites de Collatz dont les premiers termes sont des entiers impairs.

(2.4) $\qquad\qquad \mathbb{M}^T(1, \dots) = \lim_{n \to +\infty} \mathbb{M}^T(1, n)$

Le deuxième tableau contient toutes les suites de Collatz dont les premiers termes sont des entiers pairs.

(2.5) $\qquad\qquad \mathbb{M}^T(2, \dots) = \lim_{n \to +\infty} \mathbb{M}^T(2, n)$

**Exemple     2.1**

Sur la figure suivante, on représente les deux tableaux complets de Collatz d'ordre dimentionnel4 :

| 31 | 47 | 71 | 107 | 161 |
|---|---|---|---|---|
| 29 | 44 | 22 | 11 | 17 |
| 27 | 41 | 62 | 31 | 47 |
| 25 | 38 | 19 | 29 | 44 |
| 23 | 35 | 53 | 80 | 40 |
| 21 | 32 | 16 | 8 | 4 |
| 19 | 29 | 44 | 22 | 11 |
| 17 | 26 | 13 | 20 | 10 |
| 15 | 23 | 35 | 53 | 80 |
| 13 | 20 | 10 | 5 | 8 |
| 11 | 17 | 26 | 13 | 20 |
| 9 | 14 | 7 | 11 | 17 |
| 7 | 11 | 17 | 26 | 13 |
| 5 | 8 | 4 | 2 | 1 |
| 3 | 5 | 8 | 4 | 2 |
| 1 | 2 | 1 | 2 | 1 |

| 32 | 16 | 8 | 4 | 2 |
|---|---|---|---|---|
| 30 | 15 | 23 | 35 | 53 |
| 28 | 14 | 7 | 11 | 17 |
| 26 | 13 | 20 | 10 | 5 |
| 24 | 12 | 6 | 3 | 5 |
| 22 | 11 | 17 | 26 | 13 |
| 20 | 10 | 5 | 8 | 4 |
| 18 | 9 | 14 | 7 | 11 |
| 16 | 8 | 4 | 2 | 1 |
| 14 | 7 | 11 | 17 | 26 |
| 12 | 6 | 3 | 5 | 8 |
| 10 | 5 | 8 | 4 | 2 |
| 8 | 4 | 2 | 1 | 2 |
| 6 | 3 | 5 | 8 | 4 |
| 4 | 2 | 1 | 2 | 1 |
| 2 | 1 | 2 | 1 | 2 |

Figure 2: Les deux tableaux complets de Collatz d'ordre dimentionnel 4: $\mathbb{M}^T(1,4)$ et $\mathbb{M}^T(2,4)$

Pour étudier les comportements des suites infinies, on procède comme suit :

-On considère des suites finies contenues dans une même matrice complète de Collatz d'ordre dimensionnel n.

- On détermine les coefficients ou les paramètres caractéristiques d'ordre n des suites considérées,

-On fait tendre l'ordre dimensionnel vers l'infini pour obtenir des suites de longueurs infinies qui sont caractérisées par des coefficients absolus. Ces derniers correspondent aux limites des coefficients caractéristiques d'ordre n lorsque n tend vers l'infini.

**Remarque 2.1**

Noter que pour l'étude et la classification des suites infinies, on peut utiliser les suites qui s'écrivent sous la forme suivante :

(2.6) $\qquad S^c(P, n) = \left(T^1(P), T^2(P), T^2(P), \dots, T^n(P)\right)$

Dans ce cas : le premier terme est $T^1(P)$ et l'externe de la suite est $T^{n+1}(P)$.

Si la suite considérée est de longueur infini, elle est notée comme suit :

$$S^{cc}(P) = (T^1(P), T^2(P), T^2(P), \dots)$$

Ou bien on peut utiliser les suites qui s'écrivent sous la forme suivante :

$$Sy(P, n) = \left(P, T^1(P), T^2(P), T^2(P), \dots, T^n(P)\right)$$

Dans ce cas si la suite est de longueur infini on adopte la notation suivante :

$$Sy^\infty(P) = (P, T^1(P), T^2(P), T^2(P), \dots)$$

Comme on va montrer que les deux suites infinies $S^{cc}(P)$ et $Sy^\infty(P)$ possèdent les mêmes coefficients caractéristiques absolus ce qui nous permet de conclure que il n'existe pas aucune différence d'utiliser les suites ou bien les pour l'étude les comportements des suites infinies de Collatz et classification conduits aux mêmes résultats.

**2.2   Tableau structurel complet du Collatz (matrice des permutations complètes)**

C'est le tableau obtenu par remplacement de chaque suite du tableau $\mathbb{M}^T(P_0, n)$ par son vecteur de parité. En réalité le tableau contient toutes les permutations avec répétition possibles de deux entiers 0 et 1 [3]. La figure suivante correspond au tableau structurel d'ordre dimensionnel 4 noté $\mathbb{M}^s(1,4)$:

| | | | | |
|---|---|---|---|---|
| 31 | 1 | 1 | 1 | 1 |
| 29 | 0 | 0 | 1 | 1 |
| 27 | 1 | 0 | 1 | 1 |
| 25 | 0 | 1 | 1 | 0 |
| 23 | 1 | 1 | 0 | 0 |
| 21 | 0 | 0 | 0 | 0 |
| 19 | 1 | 0 | 0 | 1 |
| 17 | 0 | 1 | 0 | 0 |
| 15 | 1 | 1 | 1 | 0 |
| 13 | 0 | 0 | 1 | 0 |
| 11 | 1 | 0 | 1 | 0 |
| 9  | 0 | 1 | 1 | 1 |
| 7  | 1 | 1 | 0 | 1 |
| 5  | 0 | 0 | 0 | 1 |
| 3  | 1 | 0 | 0 | 0 |
| 1  | 0 | 1 | 2 | 0 |

Figure 3 : Tableau structurel complet d'ordre 4 de Collatz

Cette dernière propriété nous permet de classifier et de catégoriser les suites de Collatz d'un même tableau complet de Collatz (Bloc des suites) en se basant uniquement sur les vecteurs de parité relatifs aux toutes les suites contenues dans ce tableau.

**2.3   Coefficients absolus caractéristiques des suites infinies de Collatz**

On définit les coefficients caractéristiques absolus relatifs à des suites de longueurs infinies, comme les limites des coefficients principaux et secondaires d'ordre n relatifs à des suites finies lorsque les longueurs de ces suites tendent vers l'infini.

**2.3.1   Coefficient principal absolu $A_\infty(P)$**

On définit le coefficient principal absolu d'une suite infinie de Collatz de premier terme P comme suit :

(2.7) $$A_\infty(P) = \lim_{n \to +\infty} A_n(P)$$

On peut distinguer deux cas possibles pour ce coefficient comme suit :

(2.8) $$\begin{cases} \lim_{n \to +\infty} A_n(P) = \text{cte (finie)} \\ \lim_{n \to +\infty} A_n(P) = +\infty \end{cases}$$

**2.3.2   Coefficient secondaire absolu $B_\infty(P)$**

On sait qu'une suite de longueur n et de premier terme est caractérisée par un coefficient secondaire qu'on le note $B_n(P)$. Pour la limite de $B_n(P)$ lorsque n tend vers l'infini, on peut distinguer deux cas différents:

-La fonction $B(P, n)$ tend vers l'infini lorsque la longueur de la suite tend vers l'infini ceci se traduit par:

$$\lim_{n \to +\infty} B_n(P) = +\infty$$

Dans ce cas, le coefficient absolu secondaire de la suite infinie considérée s'écrit :

$$B_\infty(P) = \lim_{n \to +\infty} B_n(P) = +\infty$$

-La suite de longueur infinie et de premier terme P décrit un cycle stable de longueur k qui renferme un nombre fini des termes:

$$S^c(P, n) \to (N_1, \ldots, N_k)$$

C'est à dire que la suite considérée se comporte à l' infini comme une fonction périodique

Dans ce cas ou ne peut définir une valeur unique pour la limite de $B_n(P)$ à l' infini, en réalité cette dernière décrit aussi le même cycle $(N_1, \ldots, N_k)$ à l'infini. On peut distinguer k limites différentes :

$$N_i = \lim_{n_i \to +\infty} B_{n_i}(P) = B_{i,\infty}(P)$$

On définit le coefficient secondaire absolu $B_\infty(P)$ de la suite infinie de Collatz de premier terme P comme suit :

(2.9) $$B_\infty(P) = \min\{B_{1,\infty}, B_{2,\infty}, \ldots, B_{k,\infty}\}$$

Donc dans le deux cas, ce coefficient peut être défini comme ci-dessous :

$$B_\infty(P) = \min(\lim_{n \to +\infty} B_n(P))$$

### 2.3.3 Exterme absolu d'une suite infinie de Collatz $T^\infty(P)$

L'exterme absolu d'une suite infinie est note $T^\infty(P)$ et son expression est comme suit:

(2.10) $$T^\infty(P) = A_\infty(P) P + B_\infty(P)$$

Cette équation équivaut à :

$$T^\infty(P) = A_\infty(P) P + \min(\lim_{n \to +\infty} B_n(P))$$

**Exemple 2.2**

La suite de Collatz de premier terme 7 atteint le cycle (1,2) après un nombre fini des itérations, en effet :

$$Sy^c(7,11) = (7,11,17,26,13,20,10,5,8,4,2,1)$$

On déduit que :

$$\begin{cases} T^{10}(7) = 2 \\ T^{11}(7) = 1 \end{cases}$$

On peut déduire la limite $B_n(P)$ à l' infini qui dépend de la parité de n comme suit:

$$\begin{cases} \lim_{k \to +\infty} B_{2k}(7) = 2 \\ \lim_{k \to +\infty} B_{2k+1}(7) = 1 \end{cases}$$

Par conséquent :

$$B_\infty(7) = \min\{1,2\} = 1$$

En réalité pour toutes les suites qui atteignent 1 après un nombre finie des itérations, on peut écrire:

$$B_\infty(P) = 1$$

**Remarque        2.2**

On peut montrer que toutes les suites de Collatz qui atteignent un cycle stable possède des coefficients absolus nuls ce qui signifie que :

(2.11) $$T^\infty(P) = B_\infty(P)$$

**Corollaire        2.1**

Une suite quelconque de premier terme un entier naturel non nul P qui atteigne 1 après un nombre fini des itérations satisfaisant la condition :

(2.12) $$T^\infty(P) = 1$$

**Exemple        2.3**

On considère la suite de Collatz de premier terme 1 représentée ci-dessous sous forme d'une ligne à une infinité des cases. C'est une suite à une importance particulière puisque toutes les suites qui vérifient la conjecture de Collatz rejoignent cette suite après un certain nombre des itérations.

| 1 | 2 | 1 | 2 | 1 |   |   | 1 | 2 |   |   | 1 | 2 |   |

C'est une suite périodique constituée par une infinité de cycle (1,2).

L'expression de $T^n(1)$ dépend de la parité de n donc on peut distinguer deux expressions différentes :

$$T^n(1) = \begin{cases} T^{2k}(1) = \dfrac{3^k}{2^{2k}} + \left(1 - \left(\dfrac{3}{4}\right)^k\right) & \text{si } n = 2k \\ T^{2k-1}(1) = \dfrac{3^k}{2^{2k-1}} + 2\left(1 - \left(\dfrac{3}{4}\right)^k\right) & \text{si } n = 2k-1 \end{cases}$$

On peut distinguer deux expressions différentes pour le coefficient principal comme suit :

$$\begin{cases} A_{2k}(1) = \dfrac{3^k}{2^{2k}} \\ B_{2k-1}(1) = \dfrac{3^k}{2^{2k-1}} \end{cases}$$

De même le coefficient secondaire ou d'ajustement est définie par les deux expressions ci-dessous :

$$\begin{cases} A_{2k}(1) = 1 - \left(\dfrac{3}{4}\right)^n \\ B_{2k-1}(1) = 2\left(1 - \left(\dfrac{3}{4}\right)^n\right) \end{cases}$$

On fait tendre n vers l'infini, la suite de Collatz de longueur infinie et de premier terme 1 admet un seul coefficient absolu correspond à la limite du coefficient principal d'ordre n :

$$A_\infty(1) = \lim_{k \to +\infty}(A_{2k}(1)) = \lim_{k \to +\infty}(A_{2k-1}(1)) = 0$$

La limite du coefficient d'ajustement dépend de la parité de n comme ci-dessous :

$$\begin{cases} B_{\infty,1}(1) = \lim_{n\to+\infty}(B_{2k-1}(1)) = 2 \\ B_{\infty,2}(1) = \lim_{n\to+\infty}(B_{2k}(1)) = 1 \end{cases}$$

On écrit dans ce cas :

$$B_\infty(1) = \min\{B_{\infty,1}(1), B_{\infty,2}(1)\} = 1$$

**Corollaire 2.2**

Toutes les suites de Collatz qui vérifient la conjecture de Collatz possèdent des coefficients principaux absolus nuls et leurs coefficients secondaires absolus égaux à l'unité. ces deux relations s'écrits comme ci-dessous :

(2.13) $$\begin{cases} A_\infty = 0 \\ B_\infty = 1 \end{cases}$$

**Démonstration**

On considère une suite de premier terme un entier naturel non nul P strictement supérieur à 1 et de longueur k+n tel que :

$$T^k(P) = 1$$

C'est à dire qu'elle rejoigne la suite de premier terme 1 après k itérations successives. Cette suite est représentée sous forme d'une ligne de (k+n) cases comme suit :

| P | $T^1(P)$ | $T^2(P)$ | | | $T^{k-2}(P)$ | 2 | 1 | 2 | | | 1 | 2 | |
|---|---|---|---|---|---|---|---|---|---|---|---|---|---|

$\qquad\qquad\qquad\qquad\text{Sy}(P,k) \qquad\qquad\qquad\qquad\qquad\qquad \text{Sy}^\infty(1)$

On cherche à exprimer le coefficient $A_{n+k}(P)$ en fonction de $A_n(1)$ et à exprimer le coefficient $B_{k+n}(P)$ en fonction de $A_n(1)$ et de $B_n(1)$ donc on procède comme suit:

$$T^{n+k}(P) = T^n\left(T^k(P)\right) = T^n(1) = A_n(1) + B_n(1) = A_n(1)T^k(P) + B_n(1)$$

On sait que :

$$T^k(P) = \frac{3^{m_{k-1}(p)}}{2^k} P + B_k(P) = A_k(P)\,P + B_k(P)$$

Ce qui nous permet d'écrire:

$$T^{k+n}(P) = A_n(1)\bigl(A_k(P)\,P + B_k(P)\bigr) + B_n(1)$$
$$= A_k(P)A_n(1)P + A_n(1)B_k(P) + B_n(1)$$
$$= A_{k+n}(P)\,P + B_{k+n}(P)$$

On peut tirer les deux relations suivantes:

$$\begin{cases} A_{n+k}(P) = A_k(P)A_n(1) \\ B_{k+n}(P) = A_n(1)B_k(P) + B_n(1) \end{cases}$$

Comme $A_k(P)$ est constant alors , on écrit :

$$\lim_{n\to+\infty} A_{k+n}(P) = \lim_{n\to+\infty}\bigl(A_k(P)A_n(1)\bigr) = A_k(P)\lim_{n\to+\infty}A_n(1) = 0$$

Pour le coefficient absolu d'ajustement , on écrit :

$$\lim_{n\to+\infty} B_{k+n}(P) = \lim_{n\to+\infty}\bigl(A_n(1)B_k(P) + B_n(1)\bigr)$$
$$= B_k(P)\lim_{n\to+\infty}A_n(1) + \lim_{n\to+\infty}B_n(1)$$

or
$$B_k(P) \lim_{n \to +\infty} A_n(1) = 0$$

Il en résulte que :
$$B_\infty(P) = \lim_{n \to +\infty} B_{k+n}(P) = \lim_{n \to +\infty} B_n(1) = 1$$

**Corollaire      2.3**

On considére la suite infinie de collatz suivante de premier terme un entier naturel non nul P:

$$Sy^\infty(P) = (P_1, T^1(P), T^2(P), T^3(P), \dots)$$

Alors les coeffcients absolus de la suite considerée satisfaisant les relations suivantes :

$$\begin{cases} A_\infty(P) = A_\infty(T^1(P)) = A_\infty(T^2(P)) = \\ B_\infty(P) = B_\infty(T^1(P)) = B_\infty(T^2(P)) = \\ T^\infty(P) = T^\infty(T^1(P)) = T^\infty(T^2(P)) = \end{cases}$$

Donc on peut choisier le point de départ qui convient mieux pour déterminer les coefficients absolus d'une suite infinie.

**Proposition      2.1**

On suppose qu'elle existe une suite de Collatz de longueur infinie et de premier terme un entier naturel non nul bien déterminé P tel que son coefficient principal absolu $A_\infty(P)$ est infini alors son coefficient absolu secondaire $B_\infty(P)$ est nécessairement infini aussi.

**Démonstration**

En premier temps, noter que si on suppose qu'elle existe une suite de longueur infinie et de premier terme un entier naturel non nul bien déterminé P tel que son coefficient principal absolu $A_\infty(P)$ est infini alors cette suite doit contenir une proportion non nulle des entiers naturels impairs en particlulier elle contient au moins deux entiers impairs $P_1$ et $P_2$ (condition suffisante pour la démonstration de la proposition 2.2)

On considère la suite infinie de Collatz et on désigne par $P_1$ et $P_2$ deux termes impairs quelconques appartenant a cette suite qu'on représente sous forme dun ligne a une infinite des cases.

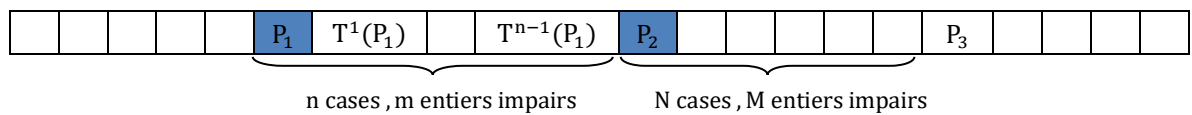

$P_2$ est l'image d'ordre n de par T, on peut écrire alors :
$$P_2 = T^n(P_1)$$

Ou encore :
$$P_2 = \frac{3^m}{2^n} P_1 + b_1$$

Comme $P_1$ est impair donc le coefficient secondaire est strictement positif.
$$b_1 > 0$$

De meme $P_3$ représente l'image d'ordre N de $P_2$, on fait référence à la ligne ci-desssus on peut écrire :

$$P_3 = \frac{3^M}{2^N} P_2 + b_2$$

Remplaçons $P_2$ par son expression dans cette dernière équation, on obtient :

$$P_3 = \frac{3^M}{2^N} P_2 + b_2 = \frac{3^M}{2^N}\left(\frac{3^m}{2^n} P_1 + b_1\right) + b_2 = \frac{3^{m+M}}{2^{N+n}} P_1 + \frac{3^M}{2^N} b_1 + b_2$$

On désigne par $P_0$ le premier terme de la suite consideree et on pose :

$$P_1 = \frac{3^{m_0}}{2^{n_0}} P_0 + b_0$$

Par conséquent:

$$P_3 = \frac{3^{m_0+m+M}}{2^{n_0+N+n}} P_0 + \frac{3^{M+m}}{2^{N+n}} b_0 + \frac{3^M}{2^N} b_1 + b_2$$

On pose :

$$b_3(N) = \frac{3^{M+m}}{2^{N+n}} b_0 + \frac{3^M}{2^N} b_1 + b_2$$

Comme la suite infinie considérée admet un coeffcient principal absolu infini, Ceci nous permet d'écrire:

$$\lim_{N \to +\infty} \frac{3^{m_0+m+M}}{2^{n_0+N+n}} = \frac{3^{m_0+m}}{2^{n_0+n}} \lim_{N \to +\infty} \frac{3^M}{2^N} = +\infty$$

Donc on déduit que :

$$\lim_{N \to +\infty} \left(\frac{3^M}{2^N}\right) = +\infty$$

On peut conclure que :

$$\lim_{N \to +\infty} b_3(N) = \lim_{N \to +\infty} \left(\frac{3^{M+m}}{2^{N+n}} b_0 + \frac{3^M}{2^N} b_1 + b_2\right) = +\infty$$

Comme on a :

$$\lim_{N \to +\infty} \frac{3^{m_0+m+M}}{2^{n_0+N+n}} = A_\infty(P_0)$$

$$\lim_{N \to +\infty} b_3(N) = B_\infty(P_0)$$

On peut conclure que si $A_\infty(P_0)$ est infini alors $B_\infty(P_0)$ est infini.

## 2.4 Généralisation sur les suites cycliques

Nous devons tout a bord faire la distinction entre les suites $\alpha-$ périodiques et les suite $\beta-$ périodiques.

**Définition    2.2**

Une suite $\alpha-$ périodique est toute suite (finie ou infinie) de Collatz qu'est caractérisée par une périodicité entière et une périodicité structurelle simultanément.

**Exemple    2.4**

La suite suivante représentée sous forme d'une ligne à 12 cases représente une suite $\alpha-$**périodique** de longueur 12. Elle est caractérisée par une périodicité numérique et une périodicité structurelle. Elle est formée par le même cycle (2,1)

| 2 | 1 | 2 | 1 | 2 | 1 | 2 | 1 | 2 | 1 | 2 | 1 |
|---|---|---|---|---|---|---|---|---|---|---|---|

**Définition 2.3**

Une suite $\beta-$périodique est toute suite (finie ou infinie) de Collatz qui possède une périodicité structurelle uniquement.

**Exemple 2.5**

La suite ci-dessous représentée sous forme d'une ligne à 12 cases représente une suite **$\beta-$périodique** de longueur 12. Elle est caractérisée par une périodicité structurelle uniquement (deux entiers impairs suivies par un entier pair) elle est formée par le même vecteur de parité (1,1,0)

| 9363 | 14045 | 21068 | 10534 | 5267 | 7901 | 11852 | 5926 | 2963 | 4445 | 6668 | 3334 |
|------|-------|-------|-------|------|------|-------|------|------|------|------|------|

**Proposition 2.2**

Soit $S^{cc}(P)$ une suite de Collatz de longueur infinie, si $S^{cc}(P)$ est une **suite $\alpha-$périodique de longueur infinie** alors elle possède nécessairement un coefficient **principal absolu** nul et un coefficient **secondaire absolu** fini égal à une valeur constante bien déterminée. Autrement :

(2.14) $$\begin{cases} A_\infty(P) = 0 \\ B_\infty(P) = cte \text{ (finie)} \end{cases}$$

**Démonstration**

On considère une suite de Collatz $\alpha-$périodique de premier terme un entier naturel non nul P et qu'on la représente comme ci-dessous sous forme d'une ligne ayant une infinité de cases.

| P |  |  | P |  |  | P |  |  | P |  |  |  |
|---|--|--|---|--|--|---|--|--|---|--|--|--|

Chaque cycle est caractérisé par un nombre des entiers impairs égal à m. et une longueur égale à n.

| P | $T^1(P)$ |  |  | $T^{n-1}(P)$ | $T^n(P) = P$ |
|---|----------|--|--|--------------|--------------|

$\underbrace{\qquad\qquad\qquad\qquad}$

n cases contiennent m entiers impairs

L'expression de $T^n(P)$ en fonction de P est comme suit :

$$T^n(P) = \frac{3^m}{2^n}P + B_n(P) = P$$

On peut écrire :

$$T^n(T^n(P)) = \frac{3^m}{2^n}\left(\frac{3^m}{2^n}P + B_n(P)\right) + B_n(P) = (\frac{3^m}{2^n})^2 P + \frac{3^m}{2^n}B_n(P) + B_n(P)$$

Par application k fois successives de $T^n$ on obtient :

$$T^n(T^n \ldots (T^n(P))) = T^{n \times k}(P) = (\frac{3^m}{2^n})^k P + B_n(P)(1 + \frac{3^m}{2^n} + \cdots + (\frac{3^m}{2^n})^{k-1})$$

La somme suivante correspond à la somme de k premiers termes d'une suite géométrique donc on peut écrire :

$$1 + \frac{3^m}{2^n} + \cdots + (\frac{3^m}{2^n})^{k-1} = \frac{1 - (\frac{3^m}{2^n})^k}{1 - \frac{3^m}{2^n}}$$

L'expression de $T^{n \times k}(P)$ d'une suite à k cycles est la suivante :

$$T^{n \times k}(P) = (\frac{3^m}{2^n})^k P + \left(\frac{1 - \left(\frac{3^m}{2^n}\right)^k}{1 - \frac{3^m}{2^n}}\right) B_n(P)$$

Comme :

$$\frac{3^m}{2^n} < 1$$

Donc on déduit que :

$$\lim_{k \to +\infty} (\frac{3^m}{2^n})^k = 0$$

De plus on a :

$$\lim_{k \to +\infty} \left(\frac{1 - \left(\frac{3^m}{2^n}\right)^k}{1 - \frac{3^m}{2^n}}\right) B_n(P) = \frac{B_n(P)}{1 - \frac{3^m}{2^n}}$$

Le coefficient absolu d'ajustement d'une suite $\alpha$ −périodique infinie de Collatz à pour expression :

$$B_\infty(P) = \frac{B_n(P)}{1 - \frac{3^m}{2^n}}$$

Comme n et m sont deux constantes, on conclut que $B_\infty(P)$ est fini.

On peut conclure que la suite $\alpha$ −périodique de longueur infinie admet un coefficient principal absolu égal à 0 et un coefficient d'ajustement absolu fini.

**Proposition     2.3**

Soit $Sy^\infty(P)$ une suite de Collatz de premier terme un entier naturel non nul P et de longueur infinie donc si $Sy^\infty(P)$ est bornée alors nécessairement elle rejoint une suite $\alpha$ −périodique si non elle tend vers l'infini.

**Démonstration**

La suite est bornée donc ils existent deux entiers naturels non nuls a et b tel que pour tout entier naturel k, on a :

$$\mathbf{a} \leq T^k(P) \leq b$$

Le nombre des entiers naturels qui comprises entre a et b est fini alors que le nombre de termes de la suite considérée est infini ce qui signifie qu'il existe au moins un terme $N_i$ de la suite $Sy^\infty(P)$ comprit entre a et b et un entier naturel non nul n tel que :

$$T^n(N_i) = N_i$$

La suite $Sy^\infty(N_i)$ est $\alpha$–périodique et les deux suites $Sy^\infty(N_i)$ et $Sy^\infty(P)$ admettent une portion infinie commune.

Dans l'autre cas il est évident que la suite est divergente et elle tend vers l'infini.

**Proposition 2.4**

Il n'exsiste pas aucun entier naturel non nul bien déterminé P tel que la suite $Sy^\infty(P)$ est une suite $\beta$–cyclique de longueur infinie.

3. **Étude de différents cas possibles et les règles générales de la classification**

La classfication des suites de collatz et la détermination des différentes proportions sont effectuées selon des règles bien déterminées.

Tableau 2 : Tableau complet fini et le tableau complet infini des collatz

| Tableau complet de collatz d'ordre dimentionnel **fini** $\mathbb{M}^T(P_0, n)$ | Tableau complet de collatz d'ordre dimentionnel **infini** $\mathbb{M}^T(P_0, \ldots)$ |
|---|---|
| n colonnes $2^n$ lignes | Infinité de colonnes et une infinité des lignes |
| Chaque suite contenue dans le tableau possède un premier terme bien déterminé | Grande proportion des suites qui font partie du tableau ne possèdent pas des premiers finis |
|  | La poroportion des suites infinies qui possèdent des premiers termes finies sont caractérisées aussi par $A_\infty =$ finie; $B_\infty$ finie; $T^\infty(P) \leq P$ (on doit prouver ces propriétés) |

Lorsqu'on fait tendre l'ordre dimensionnel vers l'infini, on obtient un tableau qui contient une infinité des suites de longueurs infinies. Les premiers termes d'une proportion de ces suites deviennent pratiquement **infinis.**

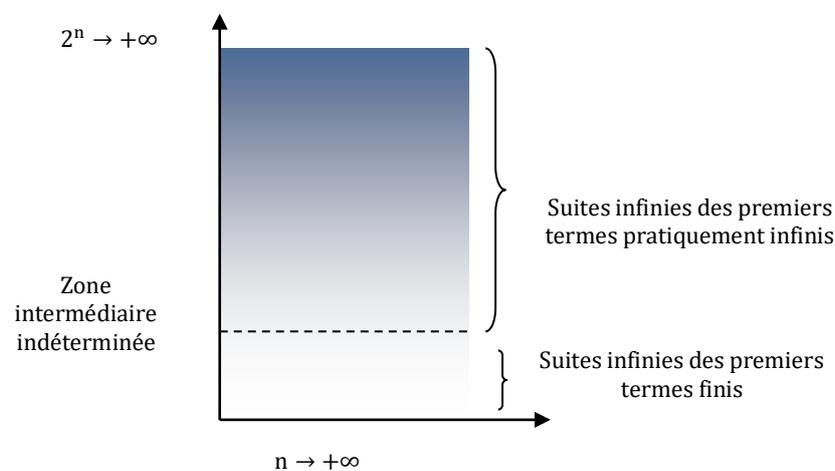

Figure 4: Comportement et répartition des suites de Collatz lorsque l'ordre dimensionnel du tableau complet (ou parfait) tend vers l'infini

On sait que pour un tableau complet de Collatz d'ordre dimensionnel infini (n → +∞), les premiers termes d'une proportion des suites contenues dans ce tableau tendent vers l'infini donc on doit faire la distinction entre les deux catégories suivantes :

-Les suites de longueurs infinies et dont les premiers termes sont des entiers finis. Elles sont qualifiées des suites R ou réalisables.

-Les suites de longueurs infinies et dont les premiers termes sont infinis ou indéterminés. Elles sont nommées des suites non réalisables NR (ou non déterminées ND).

**Exemple 3.1**

Pour tout entier naturel non nul n, on peut déterminer une suite de Collatz de longueur n et de premier terme $T^1(2^{n+1} - 1)$ contenue dans le tableau complet de Collatz d'ordre dimensionnel n comme par exemple :

$S^c(31, n = 4) = (47, 71, 107, 161)$

| 31 | 47 | 71 | 107 | 161 |
|----|----|----|-----|-----|

$S^c(127, n = 6) = (191, 287, 431, 647, 971, 1457)$

| 127 | 191 | 287 | 431 | 647 | 971 | 1457 |
|-----|-----|-----|-----|-----|-----|------|

Lorsque n tend vers l'infini, la suite $S^c(2^{n+1} - 1, n)$ n'admet pas un premier terme fini dans le tableau d'ordre dimensionnel infini.

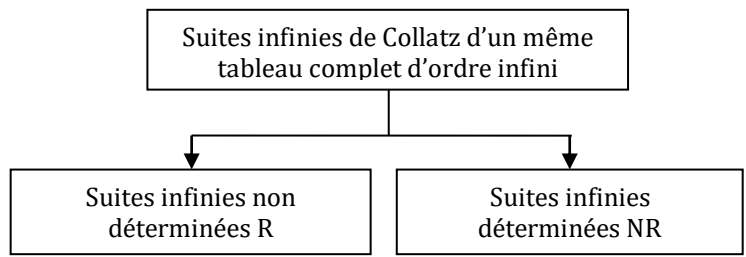

Figure 6: Les deux grandes catégories des suites infinies de Collatz

**Remarque 3.1**

Les deux types des suites réalisables et non réalisables(ou non déterminées) font partie du tableau complet de Collatz d'ordre dimensionnel infini et elles interviennent dans les calculs des différentes proportions de plus elles peuvent être classées selon les mêmes règles et les même conditions.

**Remarque 3.2**

Comme on va montrer même pour une suite infinie de premier terme infini on peut les caractériser par deux coefficients absolus un principal et un autre secondaire, on peut le qualifier de croissante ou décroissante, on se basant sur les propriétés de vecteur de parité relatif à cette suite.

**Proposition 3.1**

Toutes les suites β−périodiques de longueurs infinies sont NR, elles ne possèdent pas des termes finis.

La classification des différentes suites infinies de collatz nécessite la détermination de tous les cas possibles pour les comportements des coefficients absolus de ces suites.

### 3.1 Cas des suites réalisables (ou déterminées)

Les suites infinies de Collatz de premiers termes sont finis ou déterminées sont dites des suites **réalisables.** Elles admettent des vecteurs de parité de longueurs infinies qui sont appelés structures binaires infinies **convertibles** en suites de Collatz.

#### 3.1.1 Les suites infinies $A^+$ et les suites infinies $A^-$

On sait que les coefficients absolus principaux des suites infinies peuvent être soit finis soit infinis mais il sera plus utile parfois de classifier les suites infinies réalisables basée sur la comparaison entre le coefficient principal absolu par rapport à l'unité ce qui nous permet de distinguer deux type différent :

Les suites de type $A^+$ c'est sont les suites qui satisfaisant la condition :

(3.1) $$A_\infty(P) > 1$$

Les suites de type $A^-$ sont des suites qui satisfaisant la condition suivante:

(3.2) $$A_\infty(P) < 1$$

#### 3.1.2 Les suites infinies $B^+$ et les suites infinies $B^-$

On effectue une catégorisation des suites de Collatz basée sur le comportement des coefficients absolus d'ajustement, lorsque les suites tendent vers l'infini. On peut distinguer deux cas différents:

(1) Des suites caractérisées par des coefficients $B_\infty$ finies, on parle dans ce cas des suites $B^-$.

(2) Des suites caractérisées par des coefficients $B_\infty$ infinies (ou non majorées) , on parle dans ce cas des suites $B^+$.

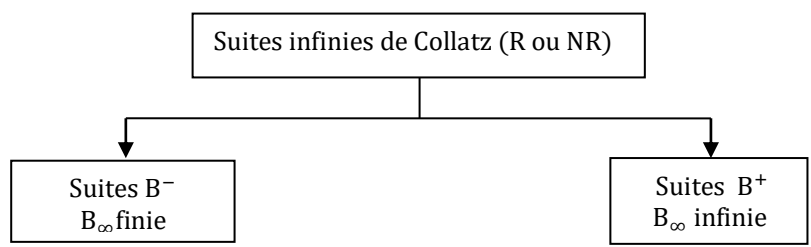

Figure 7: Classification selon la limite de la fonction $B_n$ à l'infini

#### 3.1.3 Les suites infinies $S^+$ et les suites infinies $S^-$

Une caractérisation des suites infinies de Collatz basée sur la nature de $T^\infty(P_i)$ qui peut être fini ou infini sans préciser si la suite est croissante ou décroissante.

Les suites qui vérifient $T^\infty(P_i)$ est fini sont de type $M^-$ et dans le cas contraire, on parle des suites $M^+$.

On effectue une autre catégorisation des suites de Collatz qui consiste à comparer les deux termes $T^\infty(P_i)$ et $T^1(P_i)$

Les suites $S^+$ correspondent à des suites qui vérifient :

(3.3) $$T^\infty(P) > T^1(P)$$

Les suites $S^-$ sont de suites infinies qui satisfaisant la condition :

(3.4) $$T^\infty(P) \leq T^1(P)$$

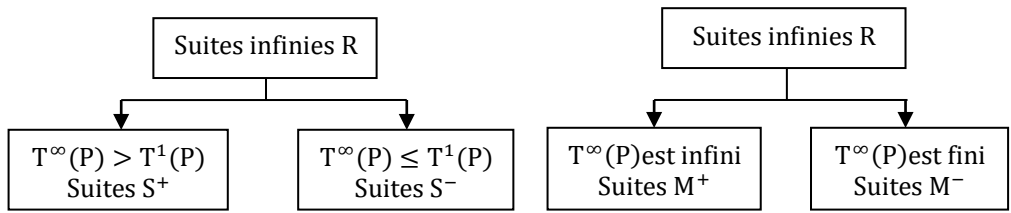

Figure 8 : classification des suites infinie réalisables basée sur le comportement de $T^\infty(P)$

**3.2 Cas des suites infinies non déterminées**

Dans la matrice parfait d'ordre dimensionnel infini, Les suites infinies de Collatz dont les premier termes sont infinis (ou indéterminés) sont dites des suites **non réalisables** (NR).
La notion des suites infinies non déterminées ou bien des structures binaires non convertibles fait appelle à la notion des séries séquentielles de Collatz qu'on va les définir dans la section suivante.

**Définition        3.1**

Une série séquentielle de Collatz est un ensemble infini des suites finies de Collatz qui remplissant un certain nombre des conditions et qu'on le note comme suit :

(3.5) $\quad Sr = (S^c(P_1, n_1), S^c(P_2, n_2), S^c(P_3, n_3), \ldots, S^c(P_k, n_k), \ldots)$

Tel que :

(3.6) $\quad S^c(P_1, n_1) \in \mathbb{M}(P_0, n_1); \; S^c(P_2, n_2) \in \mathbb{M}(P_0, n_2); \ldots; S^c(P_k, n_k) \in \mathbb{M}(P_0, n_k); \ldots$

(3.7) $$V(P_1, n_1) \subset V(P_2, n_2) \subset V(P_3, n_3) \subset \cdots \subset V(P_k, n_k) \subset \cdots$$

Avec :

- $V(P_k, n_k)$ est le vecteur de parité de la suite $Sy(P_k, n_k)$
- $\mathbb{M}(P_0, n_k)$ est la matrice parfaite de Collatz d'ordre dimensionnel $n_k$ et de générateur $P_0$ que peut être 1 ou bien 2.

**Remarque importante        3.3**

On peut utiliser les suites $Sy(P, n)$ au lieu d'utiliser des suites $S^c(P, n)$

**Définition        3.2**

(1) On définit la suite génératrice de la série séquentielle comme étant l'ensemble constituée par tous les premiers termes des suites constituants cette série:

$$\mathbb{P}(Sr) = (P_i)_{i \geq 1}$$

(2) La suite principale de la série séquentielle représente l'ensemble constitué par tous les coefficients principaux relatifs (d'ordre fini) des tous les suites de la série considérée :

$$\mathbb{A}(Sr) = (A(P_i, n_i))_{i \geq 1}$$

(3) On définit la suite secondaire de la série séquentielle comme suit :

$$\mathbb{B}(Sr) = (B(P_i, n_i))_{i \geq 1}$$

(4) On définit la suite des extermes de la série séquentielle considérée comme :

$$\mathbb{T}(Sr) = (T^{n_i+1}(P_i))_{i \geq 1}$$

La suite infinie appartenant au tableau complet d'ordre dimensionnel infini correspond à la limite suivante :

(3.8) $$S_\infty = \lim_{n_i \to +\infty} S^c(P_i, n_i)$$

La série séquentielle et sont représentées sous forme d'un tableau à une infinité des lignes comme suit :

|  | $P_i$ | $A(P_i, n_i)$ | $B(P_i, n_i)$ | $T^{n_i}(P_i)$ |
|---|---|---|---|---|
| $S^c(P_1, n_1)$ | $P_1$ | $A(P_1, n_1)$ | $B(P_1, n_1)$ | $T^{n_1}(P_i)$ |
| $S^c(P_2, n_2)$ | $P_2$ | $A(P_2, n_2)$ | $B(P_2, n_2)$ | $T^{n_2}(P_2)$ |
| $S^c(P_3, n_3)$ |  |  |  |  |
| $S^c(P_4, n_4)$ |  |  |  |  |
|  |  |  |  |  |
|  |  |  |  |  |
| $S^c(P_k, n_k)$ | $P_k$ | $A(P_k, n_k)$ | $B(P_k, n_k)$ | $T^{n_k}(P_k)$ |
|  |  |  |  |  |
|  |  |  |  |  |
|  |  |  |  |  |
| $S^{cc}(P_\infty)$ | $P_\infty$ | $A_\infty^s$ | $B_\infty^s$ | $T^\infty$ |

Figure 9: les suites caractéristiques d'une série séquentielle

Cette dernière suite est caractérisée par les coefficients absolus suivants :

Un premier terme correspond à la limite de la suite génératrice de la série séquentielle:

(3.9) $$P_\infty(S_\infty) = \lim_{n_i \to +\infty} P_i$$

Coefficient principal absolu correspond à la limite de la suite principale:

(3.10) $$A_\infty^s(S_\infty) = \lim_{n_i \to +\infty} A(P_i, n_i)$$

Coefficient secondaire absolu égal à la limite de la suite secondaire:

(3.11) $$B_\infty^s(S_\infty) = \lim_{n_i \to +\infty} B(P_i, n_i)$$

L'exterme absolu correspond à la limite suivante:

(3.12) $$T^\infty(S_\infty) = \lim_{n_i \to +\infty} T^{n_i+1}(P_i)$$

**Séries séquentielles convergentes et séries séquentielles divergentes**

Pour la limite de la suite $(P_i)_{i \geq 1}$ on peut distinguer deux cas possibles :

(1) Si $\lim_{n_i \to +\infty} P_i = P_\infty(S_\infty)$ est finie

Dans ce cas on dit que la série séquentielle considérée est convergente et la suite infinie $S_\infty$ est une suite réalisable dont le premier est un entier fini $P_\infty(S_\infty)$.

(2) Si $\lim_{n_i \to +\infty} P_i = P_\infty(S_\infty)$ est infinie

Dans ce cas on dit que la série séquentielle considérée est divergente et la suite $S_\infty$ est une suite infinie non réalisable en effet son premier terme $P_\infty$ est infini.

On s'intéresse plus particulièrement aux suites non réalisables issues des séries séquentielles divergentes.

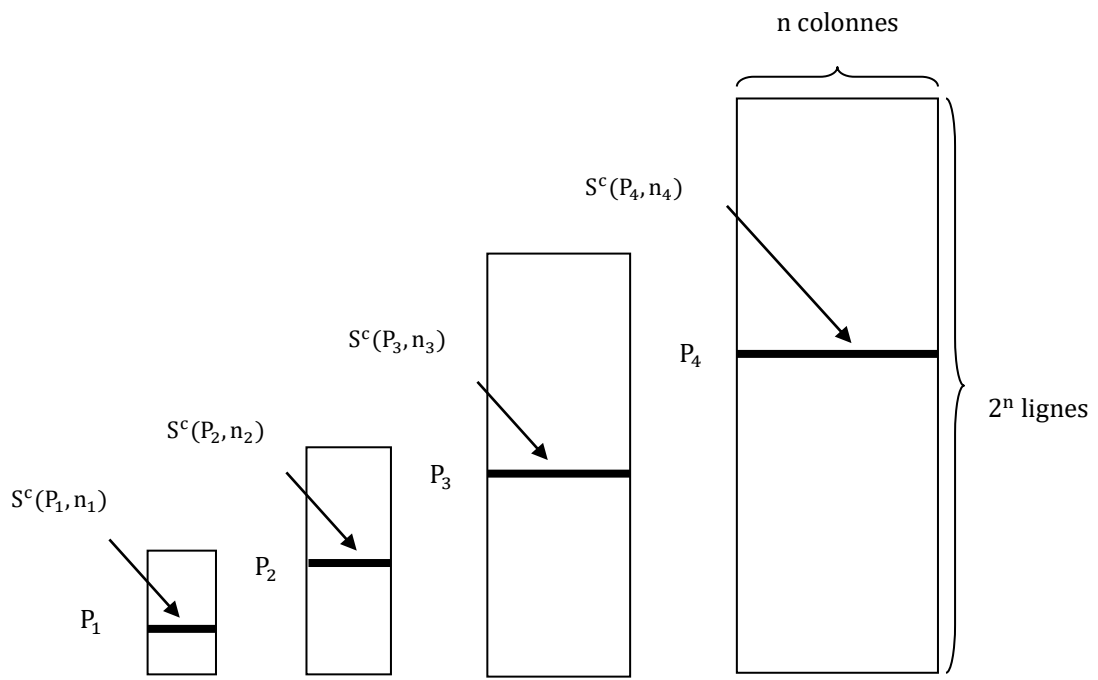

Figure 10 : série séquentielle divergente

Une suite infinie NR est dite de type $A^+$ si $A_\infty^s(S) > 1$ si non il s'agit d'une suite de type $A^-$.

Si $B_\infty^s(S)$ est fini on dit que S est de type $B^-$ si non elle s'agit d'une suite NR de type $B^+$.

**Proposition 3.2**

Toutes les suites infinies non réalisables de Collatz correspondent à une limite d'une série séquentielle divergente.

**Exemple 3.2**

Chaque tableau complet de Collatz d'ordre dimensionnel 2n, contient une suite de Collatz qui admet le vecteur de parité suivant :

$$v(n) = <\underbrace{1,1,1,1,\ldots,1}_{n\text{ fois}},\underbrace{1,0,0,0,\ldots,0,0}_{n\text{ fois}}>$$

On fait tendre n vers l'infini on obtient un vecteur de parité (ou suite binaire) de longueur infinie qu'on la note $V_\infty$

$$V_\infty = \lim_{n \to +\infty} v(n)$$

La suite binaire infinie considérée est représentée sous forme d'une ligne à une infinité de cases telles que la moitié de ces cases sont occupées par le nombre 1.

| 1 | 1 |  |  |  | 1 | 1 | 1 | 0 | 0 | 0 |  |  |  | 0 | 0 |
|---|---|---|---|---|---|---|---|---|---|---|---|---|---|---|---|

On cherche à déterminé les coefficients caractéristiques absolus caractérisant

On peut construire une suite séquentielle constituée par des suites finies $S^c(P, n)$ qui tendent vers la suite infinie considérée. Une telle série est notée, les quatre premières suites

$Sr = (S^c(3,2), S^c(23,4), S^c(15,6), S^c(351,8), ...)$

Les suites sont représentes ci-dessous avec leurs vecteurs de parités :

| 351 | 527 | 791 | 1187 | 1781 | 2672 | 1336 | 668 | 334 |   $(1,1,1,1,0,0,0,0)$   $\mathbb{M}^T(1,8)$

| 15 | 23 | 35 | 53 | 80 | 40 | 20 |   $(1,1,1,0,0,0)$   $\mathbb{M}^T(1,6)$

| 23 | 35 | 53 | 80 | 40 |   $(1,1,0,0)$   $\mathbb{M}^T(1,4)$

| 3 | 5 | 8 |   $(1,0)$   $\mathbb{M}^T(1,2)$

Figure11 : les suites des premiers rangs dans la série séquentielle considérée

Comme la suite $S^c(P_n, 2n)$ fait partie du tableau $\mathbb{M}^T(1,2n)$ donc son premier terme vérifie :
$$P_n \leq 2^{2n+1} - 1$$

Pour que les n premiers termes de cette suite soient des entiers impairs il faut qu'il existe un entier naturel impair j tel que :
$$P_n = 2^{n+1}j - 1$$

| $2^{n+1}j - 1$ | $3.2^n j - 1$ | $9.2^{n-1}j - 1$ |  |  |  | $2\times 3^n j - 1$ |

$\underbrace{\qquad\qquad\qquad\qquad\qquad\qquad\qquad\qquad}_{\text{Les n premières cases remplies des entiers impairs}}$

Ce qui nous permet de déduire :
$$P_\infty(S) = \lim_{n \to +\infty} P(n) = \lim_{n \to +\infty}(2^n j - 1) = +\infty$$

Pour une suite de longueur 2n, l'expression du coefficient principal relatif est la suivante :
$$A(P_n, 2n) = \frac{3^n}{2^{2n}}$$

Donc on peut déduire la limite suivante :
$$A_\infty^s(S) = \lim_{n \to +\infty} A(P_n, 2n) = 0$$

L'expression du coefficient secondaire relatif est la suivante :
$$B(P(n),2n) = \frac{1}{2^n}\frac{(\frac{3}{2})^n - 1}{\frac{3}{2} - 1} = \frac{1}{2^n}((\frac{3}{2})^n - 1)$$

On déduit la valeur du coefficient secondaire absolu :
$$B_\infty^s(S) = \lim_{n \to +\infty} B(P(n), 2n) = 0$$

L'expression de l'exterme de la suite de longueur 2n est comme ci-dessous :
$$T^{2n+1}(P) = \frac{1}{2^n}(3^n j - 1)$$

Donc on peut vérifier aisément que :
$$\frac{T^{2n+1}(P)}{2^{n+1} j - 1} < 1$$

Par suite on peut conclure que toutes les suites de la série considérée sont de type $S^-$ et comme $T^{2n+1}(P)$ est infinie, on déduit qu'elles s'agissent des suites $M^+$.

On peut conclure que la suite infinie qui correspond à la limite de la série séquentielle est de la catégorie $NR[A_0^- B^- S^-(M^+)]$.

**Exemple 3.3**

On peut déterminer la nature des suites $\beta$ − cycliques de longueurs infinies, on se basant uniquement à l'unité cyclique minimale qui constitue la suite infinie.

Rappeler que les suites $\beta$ − cycliques sont toutes NR. Prenons l'exemple d'une suite infinie constituée par une infinité d'un même cycle structurel qui correspond au vecteur de parité <1,1,0>. Puisque les termes de la suite sont infini donc on va la représentée par son vecteur de parité infini comme suit :

| 1 | 1 | 0 | 1 | 1 | 0 | 1 | 1 | 0 | 1 | 1 | 0 |
|---|---|---|---|---|---|---|---|---|---|---|---|

On désigne par U le vecteur de parité répétitif constitutif de cette suite :

| 1 | 1 | 0 |
|---|---|---|

Ce vecteur est caractérisé par les deux coefficients suivants :
$$\begin{cases} A^s(U) = \frac{3^2}{2^3} > 1 \\ B^s(U) = \frac{5}{8} \end{cases}$$

En exploitant les résultats trouvés lors de l'étude des suites cycliques :
$$T^{n \times k}(P) = \left(\frac{3^2}{2^3}\right)^k P + \left(\frac{1 - \left(\frac{3^2}{2^3}\right)^k}{1 - \frac{3^2}{2^3}}\right) B^s(U)$$

Le coefficient d'ajustement absolu relatif à la suite considérée est obtenu comme suit :
$$B_{3k}^s(V) = \left(\frac{1 - \left(\frac{3^2}{2^3}\right)^k}{1 - \frac{3^2}{2^3}}\right) B^s(U) \to +\infty \text{ lorsque k tend vers linfini}$$

On déduit que La suite considérée possède un coefficient secondaire absolu infini
$$B_\infty^s(V) \text{ est infini}$$

Le coefficient principal absolu de la suite considérée est le suivant :

$$A^s_\infty(V) = \lim_{k \to +\infty} (\frac{3^2}{2^3})^k = +\infty$$

La suite considérée est de la catégorie ND[$A^+B^+S^+(M^+)$]

## 4. Processus inverse de Collatz

On sait qu'une suite **finie** de Collatz peut être convertie en une suite binaire qui traduise la distribution des entiers impairs et pairs dans la suite considérée [3]. Le processus inverse consiste à considéré un vecteur de parité donné et on cherche son équivalent en suite de Collatz. Lorsqu'on fait tendre l'ordre dimensionnel d'un tableau structurel vers l'infini on obtient un tableau qui contient une infinite des vecteurs binaires des longueurs infinies.

On peut distinguer deux cas pour les vecteurs de parité des longueurs infinies dans le tableau Des vecteurs de parité possédant des rangs finis. Ces vecteurs sont appelées structure binaire infinie convertible en suites de Collatz.

| | | | | | |
|---|---|---|---|---|---|
| $2^{n+1}-1$ | | | | | |
| | | | | | |
| | | | | | |
| 31 | 1 | 1 | 1 | 1 | |
| 29 | 0 | 0 | 1 | 1 | |
| 27 | 1 | 0 | 1 | 1 | |
| 25 | 0 | 1 | 1 | 0 | |
| 23 | 1 | 1 | 0 | 0 | |
| 21 | 0 | 0 | 0 | 0 | |
| 19 | 1 | 0 | 0 | 1 | |
| 17 | 0 | 1 | 0 | 0 | |
| 15 | 1 | 1 | 1 | 0 | |
| 13 | 0 | 0 | 1 | 0 | |
| 11 | 1 | 0 | 1 | 0 | |
| 9 | 0 | 1 | 1 | 1 | |
| 7 | 1 | 1 | 0 | 1 | |
| 5 | 0 | 0 | 0 | 1 | |
| 3 | 1 | 0 | 0 | 0 | |
| 1 | 0 | 1 | 2 | 0 | |

$n \to +\infty$

Figure 12: Matrice des permutations complètes d'ordre infini

En réalité, on peut associer à un vecteur de parité quelconque deux coefficients caractéristiques qu'on les notes $A^s_n$ et $B^s_n$. Noter que :

La notation $A_n(P)$ est consacrée pour le coefficient principal d'une suite de Collatz de premier terme P et de longueur n.

La notation $A^s_n(V)$ est consacrée pour le coefficient principal d'un vecteur de parité de longueur n.

**Exemple 4.1**

On considère le vecteur de parité $v = < 1,0,0,1,1,0,1 >$ qu'on le représente sous forme d'une ligne à 7 cases comme suit:

| 1 | 0 | 0 | 1 | 1 | 0 | 1 |

Ce vecteur est caractérisé par les deux coefficients suivants :

$$A_n^s(v) = \frac{81}{128}, B_n^s(v) = \frac{211}{128}$$

Des vecteurs qui possèdent des rangs pratiquement infinis. Cette deuxième catégorie des vecteurs sont appelées structures binaires infinies non convertibles en suites de Collatz.

**Définition 4.1**

Une structure binaire infinie b est dite convertible en suite de Collatz si elle existe une suite de Collatz S de longueur infinie et de premier terme un entier naturel non nul bien définie P tel que le vecteur de parité de S est identique à la structure considérée b.

**Exemples 4.2**

L'exemple suivant correspond à une structure infinie convertible qui possède les coefficients absolus suivants :

$A_\infty = 0, B_{1,\infty} = 1; B_{2,\infty} = 2$

| 1 | 1 | 1 | 0 | 1 | 0 | 0 | 1 | 0 | 0 | 0 | 1 | 0 | 1 | 0 | 1 |

Elle s'agit d'une structure convertible en suite de Collatz en effet elle correspond à la suite infinie $Sy^{cc}(7)$ suivante :

| 7 | 11 | 17 | 26 | 13 | 20 | 10 | 5 | 8 | 4 | 2 | 1 | 2 | 1 | | |

**Définition 4.2**

On définit une structure binaire non convertible en suite de Collatz comme une suite binaire de longueur infinie tel qu'il n'existe pas aucun entier naturel non nul qui peut générer une suite de Collatz dont le vecteur de parité est identique à la structure considérée.

**Exemple 4.3**

On donne ici deux exemples pour les vecteurs non convertibles en suites de Collatz.

(1) le premier cas correspond à une structure cyclique infinie constituée par la même séquence 011 répétée une infinité des fois comme suit :

$$(011011011011011 \dots)$$

On peut la représenter sous forme d'une ligne chromatique :

Cette suite binaire infinie fait partie du tableau complet d'ordre dimensionnel infini et elle présente les propriétés suivantes :

$A_\infty^s > 1$(infini), $B_\infty^s$ infini, $P_\infty$ premier terme infini donc elle correspond à une structure binaire infinie non convertible en suites de Collatz. Plus précisément, elle s'agit d'un de type

$$NC(S^+A^+B^+)$$

(2) Le dernier exemple est l'exemple le plus simple, il correspond à la suite de premier terme $2^{n+1} - 1$ dans un tableau complet de Collatz d'ordre dimensionnel n. Le vecteur de parité est constitué uniquement du nombre 1 répété n fois. Lorsqu'on fait tendre n vers l'infini, le premier terme de la suite considérée dans le tableau d'ordre dimensionnel infini devient pratiquement infini. On peut montrer aisément les propriétés suivantes :

$A_\infty^s > 1$(infini), $B_\infty^s$ infini, $P_\infty$ infini , elle s'agit d'une structure binaire infinie non convertible de type $NC(S^+A^+B^+)$

**Remarque 4.1**

Les deux suites infinies qui font partie du tableau complet infini et qui correspondent à ces deux structures binaires sont deux suites non réalisables possédant les mêmes coefficients absolus.

**Remarque 4.2**

Si par exemple deux ou plusieurs conditions ne peuvent pas être remplies simultanément par aucune suite infinie de Collatz, la catégorie qui correspond à l'ensemble des suites qui doivent vérifier ces conditions est un ensemble vide.

## 5. Classification des suites infinies de Collatz

L'organigramme suivant montre les différents types des suites infinies réalisables de Collatz :

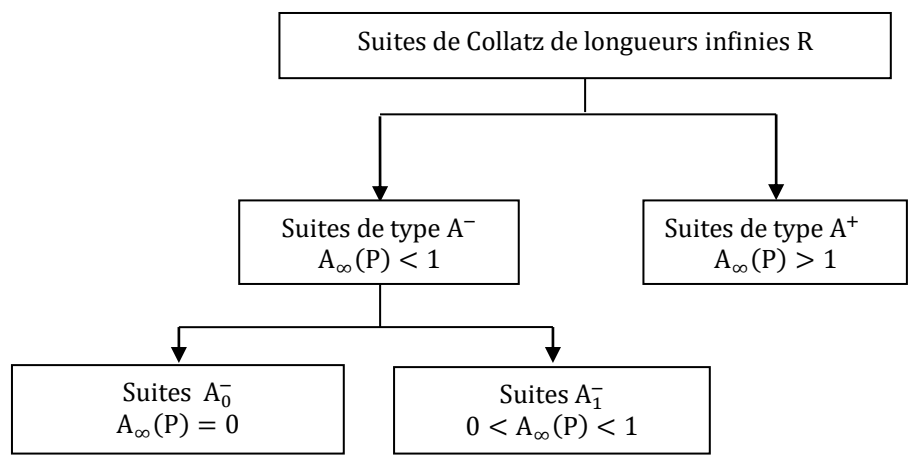

Figure 13: Différents types des suites infinies réalisables de Collatz

Puisque les suites non réalisables ne possèdent pas des premiers finis donc on fait classer chaque type de suites (R et NR) l'un indépendamment de l'autre:

### 5.2.1 Suites infinies réalisables

Elles correspondent à la proportion des suites infinies du tableau complet d'ordre infini dont leurs premiers termes sont des entiers naturels bien définies.

Tableau 3 : Classification 1 des suites **infinies réalisables** de Collatz (en gris foncée cas impossibles)

| | Suites infinies de type R | | | | | |
|---|---|---|---|---|---|---|
| | $A_\infty < 1$ Catégorie $A^-$ | | | | $A_\infty > 1$ Catégorie $A^+$ | |
| | $A_\infty = 0$ Catégorie $A_0^-$ | | $0 < A_\infty < 1$ Catégorie $A_1^-$ | | | |
| | $B_\infty$ fini | $B_\infty$ infini | $B_\infty$ fini | $B_\infty$ infini | $B_\infty$ fini | $B_\infty$ infini |
| $T^\infty(P) \leq T^1(P)$ | $RS^- A_0^- B^-$ | $RS^- A_0^- B^+$ | $R\,S^- A_1^- B^-$ | $R\,S^- A_1^- B^+$ | $RS^- A^+ B^-$ | $R\,S^- A^+ B^+$ |
| $T^\infty(P) > T^1(P)$ | $R\,S^+ A_0^- B^-$ | $RS^+ A_0^- B^+$ | $RS^+ A_1^- B^-$ | $RS^+ A_1^- B^+$ | $R\,S^+ A^+ B^-$ | $R\,S^+ A^+ B^+$ |

Les signification de différentes colorations des cases sont présentées ci-dessous :

☐ Une catégorie qu'on peut vérifier qu'elle est non vide.

☐ Des catégories qu'on croit qu'elles sont vides mais non pas encore prouvées.

☐ Elles correspondent à des catégories vides et qui sont prouvées qu'elles sont vides

Par conclusion, pour les suites infinies réalisables, on peut distinguer 12 catégories différentes dont 7 catégories sont des ensembles vides c'est à dire qu'elle n'existe pas aucune suite de premier terme bien définie que peut appartenir à l'une de ces catégories.

Au lieu de faire comparer $T^\infty(P)$ et $P$, on peut utiliser les deux condition suivantes :

$T^\infty(P)$ fini ($M^-$) ou bien $T^\infty(P)$ est infini ($M^+$)

Tableau 4 : Classification 2 des suites **infinies réalisables** de Collatz

| | Suites infinies de type R | | | | | |
|---|---|---|---|---|---|---|
| | $A_\infty < 1$ Catégorie $A^-$ | | | | $A_\infty > 1$ Catégorie $A^+$ | |
| | $A_\infty = 0$ Catégorie $A_0^-$ | | $0 < A_\infty < 1$ Catégorie $A_1^-$ | | | |
| | $B_\infty$ fini | $B_\infty$ infini | $B_\infty$ fini | $B_\infty$ infini | $B_\infty$ fini | $B_\infty$ infini |
| $T^\infty(P)$ finie | $RM^- A_0^- B^-$ | $RM^- A_0^- B^+$ | $R\,M^- A_1^- B^-$ | $R\,M^- A_1^- B^+$ | $RM^- A^+ B^-$ | $R\,M^- A^+ B^+$ |
| $T^\infty(P)$ infinie | $R\,M^+ A_0^- B^-$ | $RM^+ A_0^- B^+$ | $RM^+ A_1^- B^-$ | $RM^+ A_1^- B^+$ | $R\,M^+ A^+ B^-$ | $R\,M^+ A^+ B^+$ |

Dans ce tableau on a 8 categories vides, la categorie supplementaire correspond $R\,[M^+ A_0^- B^-]$. Une suite infinie R de collatz ne peut pas être $M^+ A_0^- B^-$ en effet comme on a :

$$T^\infty(P) = A(P, n)P + B(P, n)$$

Donc si $A(P, n)$ et $B(P, n)$ sont finies alors nécessairement $T^\infty(P)$ est finie aussi.

**Remarque importante    5.1**

Pour démontrer la conjecture de Collatz, on doit montrer que toutes les suites infinies réalisables (ou que toutes les structures binaires infinies convertibles) sont toutes appartenant à la catégorie $RS^-A_0^-B^-$ de plus on doit démontrer l'unicité de cycle (1,2). Autrement, on doit démontrer que les 4 catégories colorées en gris clair correspondent à des ensembles vides. En réalité l'unité de cycle (1,2) signifie l'impossibilité de l'existence d'aucune suite de type $RS^+A_0^-B^-$ ce qui nous permet de conclure si les conditions suivantes sont remplies :

(5.1) $\begin{cases} \{RS^+A_0^-B^+\} = \emptyset \\ \{RS^+A_1^-B^+\} = \emptyset \\ \{RS^+A^+B^+\} = \emptyset \\ \text{Le cycle (1,2) est le seul cycle stable par la transformation de collatz} \end{cases}$

Donc la conjecture de Collatz est vérifiée par n'importe quelle suite de premier terme un entier naturel non nul bien déterminé.

**Proposition    5.1**

Les ensembles qui correspondent aux catégories suivantes sont des ensembles vides, elles ne contiennent aucune suite de Collatz de longueur infinie et de premier terme un entier naturel non nul bien définie :

$R[S^-A_0^-B^+], R[S^-A_1^-B^-], R[S^+A_1^-B^-], R[S^-A_1^-B^+], R[S^-A^+B^-], R[S^-A^+B^+]$ et $R[S^+A^+B^-]$

**Démonstration**

*Impossibilité de l'existence des suites infinies $RA_1^-S^-B^+$ et de $RA_0^-S^-B^+$*

Prenons le cas des suites caractérisées par :

$$0 < A_\infty(P) < 1 \; et \; P \; est \; fini$$

Elles s'agissent des suites infinies de la catégorie $[RA_1^-]$

On sait que :

$T^\infty(P) = A_\infty(P)P + B_\infty(P)$

Si $B_\infty(P)$ est infini alors $T^\infty(P)$ est infini ce qui nous permet de conclure que $T^\infty(P) > P$ et par conséquent la suite considérée est de type $S^+$

$$S^c(P) \in [RA_1^-B^+] \Rightarrow S^c(P) \in [RA_1^-S^+B^+]$$

Ce qui signifie qu'elle n'existe pas aucune suite de type $RA_1^-S^-B^+$

De même si $A_\infty(P) = 0$ et P est fini la suite considérée est de type $[RA_1^-]$. Si $B_\infty(P)$ est infini alors $T^\infty(P)$ est infini ce qui nous permet de conclure que $T^\infty(P) > P$ et par conséquent la suite considérée est de type $S^+$

$$S^c(P) \in [RA_0^-B^+] \Rightarrow S^c(P) \in [RA_0^-S^+B^+]$$

Ce qui signifie qu'elle n'existe pas aucune suite de type $RA_0^-S^-B^+$ (cas impossible).

*Impossibilité de l'existence des suites infinies $RS^-A_1^-B^-$ et de $RS^+A_1^-B^-$*

Si on suppose qu'elle existe une suite $Sy^{cc}(P)$ de longueur infinie, de type $\mathbf{A_1^-}$ et de premier terme un entier naturel P bien défini donc elles s'agissent des suites infinies de la catégorie $[RA_1^-]$

On sait que :

$$T^\infty(P) = A_\infty(P)P + B_\infty(P)$$

Si $B_\infty(P)$ est fini donc nécessairement $T^\infty(P)$ est finie et par conséquent la suite considérée est bornée d'après le théorème 2.3, la suite considérée rejoint nécessairement une suite $\alpha$ – périodique et par conséquent le coefficient principal absolu de la suite $Sy^{cc}(P)$ est nul absurde puisque cette suite ne possède pas un coefficient principal absolu nul. On peut conclure qu'elle n'existe pas aucune suite de premier terme un entier bien définie qui peut être à la fois $A_1^-$ et $B^-$ et par conséquent $[RA_1^-B^-]$ est vide et elle n'existe pas aucune suite infinie de type $RA_1^-B^-S^+$ ou bien $RA_1^-B^-S^-$.

***Impossibilité de l'existence des suites infinies $RS^-A^+B^-$ et de $RS^-A^+B^-$***

Concernant ces deux catégories, il est évident que si une suite appartient à la catégorie $[A^+]$ donc elle appartient nécessairement à la catégorie $[S^+]$ en effet :

$$T^\infty(P) = A_\infty T^1(P) + B_\infty(P)$$

Donc on peut écrire :

$$\text{si } A_\infty > 1 \Rightarrow T^\infty(P) > T^1(P)$$

***Impossibilité de l'existence des suites infinies $RS^+A^+B^-$***

$$T^\infty(P) = A_\infty P + B_\infty$$

Avec les conditions suivantes:

$$\begin{cases} A_\infty > 1 \\ P \text{ est fini} \\ B_\infty \text{ est fini} \end{cases}$$

Si $T^\infty(P)$ est fini donc nécessairement la suite $S^{cc}(P)$ rejoint une suite $\alpha$ – périodique et par conséquent le coefficient principal absolu $A_\infty$ est pratiquement nul ce qui est absurde puisque $A_\infty$ est non nul donc $T^\infty(P)$ ne peut pas être fini, il doit être infini est par conséquent $A_\infty$ est infini aussi. En utilisant la proposition 2.1 qui affirme que toutes suites infinies admettant des coefficients absolus principaux infinis, elles admettent aussi nécessairement des coefficients absolus secondaires infinis. Ce qui nous permet de conclure que si $A_\infty$ est infini alors $B_\infty$ est infini aussi ceci prouve l'impossibilité de l'existence des suites $RS^+A^+B^-$.

### 5.2.2 Suites infinies non réalisables NR (ou non déterminées ND)

Elles correspondent à la proportion des suites de Collatz de longueurs infinies contenues dans le tableau complet de Collatz d'ordre dimensionnel infini et qui ne possèdent pas des termes finis ou déterminés.

Tableau 5: Classification des suites **infinies non réalisables** de Collatz

| Suites infinies de type NR ou ND (suites de premiers termes infinis ou indéterminés) | | | | | |
|---|---|---|---|---|---|
| $A_\infty^s < 1$ | | | | $A_\infty^s > 1$ | |
| $A_\infty^s = 0$ | | $0 < A_\infty^s < 1$ | | | |
| $B_\infty^s$ fini | $B_\infty^s$ infini | $B_\infty^s$ fini | $B_\infty^s$ infini | $B_\infty^s$ fini | $B_\infty^s$ infini |
| $ND(A_0^- B^-)$ | $ND(A_0^- B^+)$ | $ND(A_1^- B^-)$ | $ND(A_1^- B^+)$ | $ND(A^+ B^-)$ | $ND(A^+ B^+)$ |

La classification de ce genre des suites est moins importante que la classification des suites réalisables.

### 5.3 Classification structurelle relative aux suites infinies de Collatz

Comme chaque suite binaire infinie correspond à une suite infinie réalisable ou non réalisable donc la classification des structures binaires infinies repose sur les mêmes principes. A titre d'exemple si on désigne par v un vecteur de parité de longueur infinie alors :

-Si $A_\infty^s(v) > 1$ on dit que v est de type $A^+$ si non, on dit qu'il est de type $A^-$.

-Si $B_\infty^s(v)$ est finie on dit que v est de type $B^-$ si non on dit qu'il est de type $B^+$.

La première classification faire la distinction entre les structures binaires non convertibles et celles qui sont convertibles (classification similaire aux cas des suites infinies réalisables et suites infinies non réalisables) :

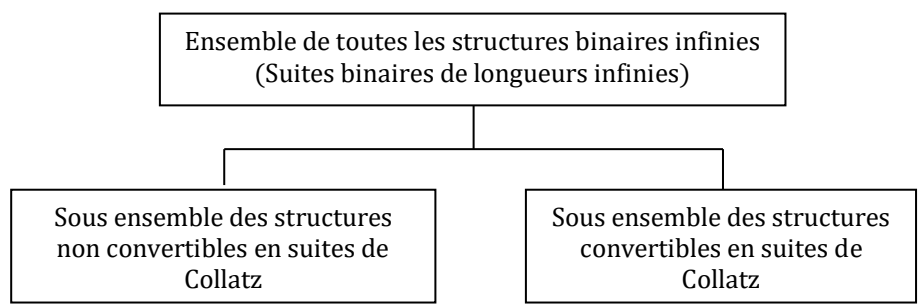

Figure 14: Classification des structures binaires infinies

### 5.3.1 Classification des structures binaires infinies convertibles

La classification des structures binaires infinies convertibles est présentée sur le tableau suivant :

Tableau 6 : Classification des structures **infinies convertibles**

|  | Structures infinies convertibles ||||||
|---|---|---|---|---|---|---|
|  | $A_\infty < 1$ |||| $A_\infty > 1$ ||
|  | $A_\infty = 0$ || $0 < A_\infty < 1$ || ||
|  | $B_\infty$ fini | $B_\infty$ infini | $B_\infty$ fini | $B_\infty$ infini | $B_\infty$ fini | $B_\infty$ infini |
| $T^\infty(P)$ $\leq T^1(P)$ | $C(S^-A_0^-B^-)$ | $C(S^-A_0^-B^+)$ | $C(S^-A_1^-B^-)$ | $C(S^-A_1^-B^+)$ | $C(S^-A^+B^-)$ | $C(S^-A^+B^+)$ |
| $T^\infty(P)$ $> T^1(P)$ | $C(S^+A_0^-B^-)$ | $C(S^+A_0^-B^+)$ | $C(S^+A_1^-B^-)$ | $C(S^+A_1^-B^+)$ | $C(S^+A^+B^-)$ | $C(S^+A^+B^+)$ |

Ces structures correspondent à tous les vecteurs de parité relatives à toutes les suites de Collatz qui possèdent des premiers termes bien définis donc elles possèdent la même classification et les mêmes propriétés de suites infinies réalisables présentées dans le tableau2. Les mêmes démonstrations pour les cas impossibles des suites réalisables restent valables pour les structures binaires infinies convertibles.

**5.3.2 Classification des structures binaires infinies non convertibles**

Pour les structures binaires infinies non convertibles, elles correspondent aux vecteurs de parité de la proportion des suites de Collatz dont leurs premiers termes tendent vers l'infini lorsque l'ordre dimensionnel du tableau complet de Collatz tend vers l'infini. Le tableau ci-dessous montre les différentes catégories de ce genre des structures binaires infinies.

Tableau 7 : Différentes catégories des structures **infinies non convertibles**

|  | Structures binaires infinies non convertibles ||||||
|---|---|---|---|---|---|---|
|  | $A_\infty^s < 1$ |||| $A_\infty^s > 1$ ||
|  | $A_\infty^s = 0$ || $0 < A_\infty^s < 1$ || ||
|  | $B_\infty^s$ fini | $B_\infty^s$ infini | $B_\infty^s$ fini | $B_\infty^s$ infini | $B_\infty^s$ fini | $B_\infty^s$ infini |
| $T^\infty(P_\infty)$ $\leq T^1(P_\infty)$ | $NC(S^-A_0^-B^-)$ | $NC(S^-A_0^-B^+)$ | $NC(S^-A_1^-B^-)$ | $NC(S^-A_1^-B^+)$ | $NC(S^-A^+B^-)$ | $NC(S^-A^+B^+)$ |
| $T^\infty(P_\infty)$ $> T^1(P_\infty)$ | $NC(S^+A_0^-B^-)$ | $NC(S^+A_0^-B^+)$ | $NC(S^+A_1^-B^-)$ | $NC(S^+A_1^-B^+)$ | $NC(S^+A^+B^-)$ | $NC(S^+A^+B^+)$ |

Ces structures ne sont pas convertibles ou plutôt elles correspondent à des suites dont les premiers termes sont infinis lorsque les longueurs des suites tendent vers l'infini.

**Notations      5.1**

On désigne par : plutôt

-$[A_0^-]$ l'ensemble des suites de type $A_0^-$.

-$[A_1^-]$ l'ensemble des suites de type $A_1^-$.

-$[S^+]$ l'ensemble qui renferme toutes les suites infinies de Collatz $S^+$

-[$S^-$] L'ensemble qui renferme toutes les suites infinies de Collatz $S^-$

Dans tous ce qui suit, on va traiter quelques exemples simples pour la détermination et la qualification des structures infinies.

**Exemple 5.1**

On définit une structure infinie qu'on la note V caractérisée par deux parties la première est constituée par une infinité des cases blanches et une partie est constituée par une infinité des cases bleus. On construit une série d'une infinité des séquences finies tel que chaque séquence est incluse dans l'autre et l'ensemble des séquences tendent vers la structure infinie considérée.

$$(0,0,0,\ldots,0,0,1,1,1,\ldots,1,1)$$

Théoriquement, cette structure binaire infinie peut être représentée sous forme d'une ligne a une infinité des cas dont la moitie sont des zéro (occupes les premiers rangs) et l'autre moitie est occupée par le nombre 1.

$\infty \leftarrow$ | 0 | 0 | | 0 | 0 | 0 | 1 | 1 | | 1 | 1 | $\to \infty$

Pour déterminer les coefficients absolus de la structure binaire infinie considérée, on peut construire une série constituée des suites $S^c(P,n)$ (au lieu des suites $Sy(P,n)$) extraites de la matrice parfaite d'ordre n

$$N(n) = \begin{cases} \frac{1}{3}(2^{2n+1} - 2^{n+1} - 1) & \text{si n est impair} \\ \frac{1}{3}(3 \times 2^{2n+1} - 2^{n+1} - 1) & \text{si n est pair} \end{cases}$$

Le premier terme d'une séquence de longueur 2n est donnée par :

$$T^1(N) = \begin{cases} 2^{2n} - 2^n & \text{si n est impair} \\ 3 \times 2^{2n} - 2^n & \text{si n est pair} \end{cases}$$

Pour n=5

$N(5)$  $T^1(N)$

| 661 | 992 | 496 | 248 | 124 | 62 | 31 | 47 | 71 | 107 | 161 |

Par exemple pour n=1,2,3 et 4, les quatre premières séquences sont les suivantes :

| 1 | 2 | 1 |

| 29 | 44 | 22 | 11 | 17 |

| 37 | 56 | 28 | 14 | 7 | 11 | 17 |

| 501 | 752 | 376 | 188 | 94 | 47 | 71 | 107 | 161 |

Figure 15: Les premières suites de la série séquentielle de la suite infinie considérée

Les expressions de $T^{2n+1}(N(n))$ sont les suivantes :

$$T^{2n+1}(N(n)) = \begin{cases} 3^n - 1 & \text{si n est impair} \\ 3^{n+1} - 1 & \text{si n est pair} \end{cases}$$

Noter que $T^{2n+1}(N)$ ne fait pas partie de la séquence considérée.

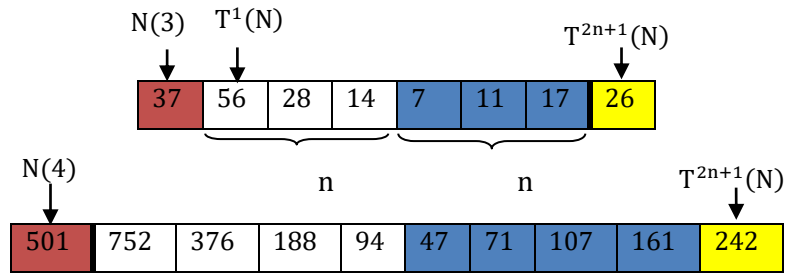

L'expression générale de $T^{2n+1}(N)$ en fonction de N est comme suit :

$$T^{2n+1}(N) = \frac{3^{n+1}}{2^{2n+1}} N(n) + \left(\frac{3^n}{2^{2n+1}} + \frac{3^n}{2^n} - 1\right)$$

On peut tirer les expressions de deux coefficients caractéristiques comme suit :

$$\begin{cases} A_n = \dfrac{3^{n+1}}{2^{2n+1}} \\ B_n = \dfrac{3^n}{2^{2n+1}} + \dfrac{3^n}{2^n} - 1 \end{cases}$$

Les coefficients absolus correspondent à la structure V de longueur infinie sont les suivants:

$$A_\infty = \lim_{n \to +\infty} \left(\frac{3^{n+1}}{2^{2n+1}}\right) = 0$$

$$B_\infty = \lim_{n \to +\infty} \left(\frac{3^n}{2^{2n+1}} + \frac{3^n}{2^n} - 1\right) = +\infty$$

Quelque soit n la relation suivante est vérifiée :

$$T^{2n+1}(N(n)) < T^1(N(n))$$

On peut conclure que la structure infinie considérée possède :

-Un coefficient principal absolu nul.

-Un coefficient d'ajustement absolu infini.

-Les séquences sont toutes de type $S^-$.

Elle est non convertible puisque les premiers termes des séquences finies tendent vers l'infini lorsque les longueurs de ces séquences tendent vers l'infini. On peut conclure qu'elle s'agit d'une structure $NC(A_0^- B^+ S^-)$.

**Exemple        5.2**

On considère la structure infinie suivante qu'on la note $S_\infty$:

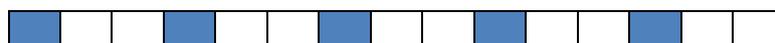

C'est une structure infinie cyclique constituée par une infinité de l'unité séquentielle cyclique suivante qu'on la note U:

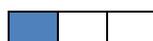

Cette séquence unitaire peut être caractérisée par les deux coefficients suivants :

$$A^s(U) = \frac{3^1}{2^3}, B^s(U) = \frac{1}{8}$$

Ceci nous permet d'écrire :

$$T^4(P) = \frac{3}{8} T^1(P) + \frac{1}{8}$$

On peut travailler avec les suites $S^c(P)$ (et non $Sy(P,n)$) extraites de la matrice complète $\mathbb{M}^T(1, n)$

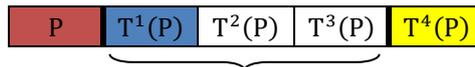

Pour une séquence composée de k unités cycliques :

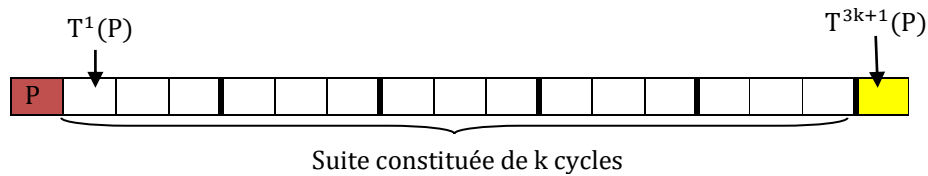

Suite constituée de k cycles

$$T^{3k+1}(P) = \left(\frac{3}{8}\right)^k T^1(P) + \frac{1}{8}\left(1 + \frac{3}{8} + \cdots \left(\frac{3}{8}\right)^{k-1}\right)$$

L'expression du coefficient relatif $B_k^s$ (d'ordre k) est la suivante :

$$B_{3k}^s(V) = \frac{1}{8}\left(1 + \frac{3}{8} + \cdots \left(\frac{3}{8}\right)^{k-1}\right) = \frac{1 - \left(\frac{3}{8}\right)^k}{1 - \frac{3}{8}} = \frac{1}{5}\left(1 - \left(\frac{3}{8}\right)^k\right)$$

Pour tout entier naturel non nul k, on a :

$$T^{3k+1}(P) = \left(\frac{3}{8}\right)^k T^1(P) + \frac{1}{5}\left(1 - \left(\frac{3}{8}\right)^k\right) < \left(\frac{3}{8}\right)^k T^1(P) + \frac{1}{5} < T^1(P)$$

On peut conclure que toutes les séquences finies de la structure considérée sont de type $S^-$

Dans ce cas on dit que la structure infinie considérée est de type $S^-$

On fait tendre k vers l'infini, on obtient :

$$\lim_{k \to +\infty} \left(\frac{1 - \left(\frac{3}{8}\right)^k}{1 - \frac{3}{8}}\right) = \frac{8}{5} \Rightarrow B_\infty^s(V) = \frac{1}{5}$$

De plus on a :

$$\lim_{k \to +\infty} \left(\frac{3}{8}\right)^k = 0$$

On déduit que :

$$\begin{cases} A_\infty^s(V) = 0 \\ B_\infty^s(V) = \frac{1}{5} \end{cases}$$

La structure infinie considérée est de la catégorie NC[$A_0^- B^- S^-$]

Les suites ci-dessous sont quelques séquences de différentes longueurs appartenant à la série constructive de la structure infinie considérée.

| 3 | 5 | 8 | 4 | | | | | | |
|---|---|---|---|---|---|---|---|---|---|
| 51 | 77 | 116 | 58 | 29 | 44 | 22 | | | |
| 819 | 1229 | 1844 | 922 | 461 | 692 | 346 | 173 | 260 | 130 |

Figure 16: Quelques suites de la série séquentielle de la suite cyclique infinie $s_\infty$

La structure binaire infinie V considérée est non convertible en suite de Collatz. Dans le tableau parfait de Collatz d'ordre dimensionnel infini, cette structure n'admette pas un rang fini, son rang tend vers l'infini. Du ce fait on dit qu'elle s'agit d'une structure binaire infinie non convertible.

**Remarque importante     5.2**

Lorsqu'on parle d'une catégorie donnée (catégorie [$A^+$] ou bien [$A^-$] ou bien un autre type parmi les types déjà définis) des suites de Collatz de longueurs infinies donc cette catégorie renferme toutes les structures convertibles et non convertibles de longueurs infinies et de même type $A^+$ ou $B^-$ ou n'importe quel autre type. De même lorsqu'on parle par exemple de la proportion des suites cette proportion est calculée par rapport au nombre total de toutes les structures C et NC du tableau complet de Collatz d'ordre dimensionnel fini ou infini.

## 6. Détermination des Différentes proportions

Dans tous ce qui suit on adopte les notations suivantes :

$r_a(S^-)$ La proportion absolue des suites de longueurs infinies de type $S^-$.

$r_a(S^+)$ La proportion absolue des suites de longueurs infinies de type $S^+$.

$r_a(A^+)$ La proportion absolue des suites de longueurs infinies de type $A^+$.

$r_a(A^-)$ La proportion absolue des suites de longueurs infinies de type $A^-$.

$r_a(R)$ La proportion absolue des suites de longueurs infinies réalisables.

$r_a(NR)$ La proportion absolue des suites de longueurs infinies non réalisables.

### 6.1 Proportion des suites infinies réalisables

D'après le théorème 8.2 [4] lorsque l'ordre dimensionnel n du tableau complet de Collatz n tend vers l'infini la proportion de suites R tend vers 0 et celle de suite NR tend vers 1.

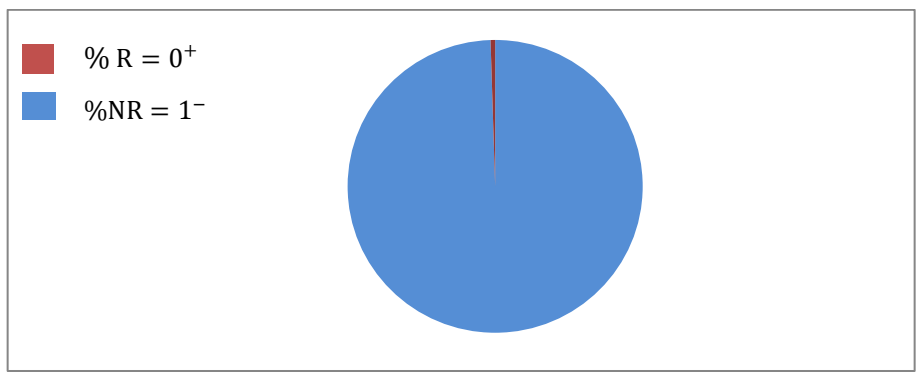

Figure 17: Répartition des proportions absolues des suites infinies de Collatz

La quasi-totalité des suites de Collatz de longueurs infinies sont non réalisables ou non déterminées c'est à dire qu'elles possèdent des termes pratiquement infinis elles admettent des vecteurs de parité de longueurs infinies non convertibles.

### 6.2 Proportion des suites infinies de type $A^-$

D'après le théorème 8.2 [4] lorsque l'ordre dimensionnel n du tableau complet de Collatz n tend vers l'infini la proportion de suites $A^-$ tend vers 1 et celle de suite $A^+$ tend vers 0. Une démonstration détaillée basée sur ces tableaux d'arrangement binaires complets nous permet d'exprimer le nombre des suites $A^+$ ou bien $A^-$ en fonction de n puis on montre que la proportion de suite $A^+$ tend vers 0 lorsque n tend vers l'infini et par conséquent la proportion des suites de type $A^-$ tend vers 1.

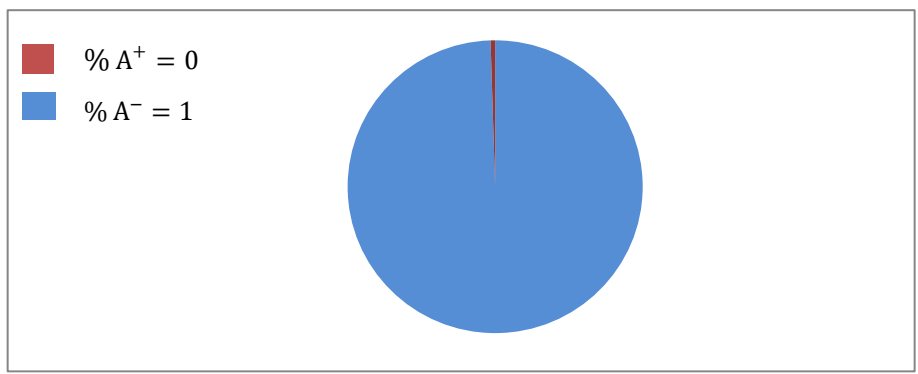

Figure 18: Répartition des proportions absolues des suites infinies de Collatz

La quasi-totalité des suites de Collatz de longueurs infinies sont de type $A^-$. C'est à dire qu'elles possèdent des coefficients principaux absolus strictement inférieurs à 1.

### 6.3  Proportion des suites infinies de type $S^-$

D'après le théorème 7.1 [4] lorsque n tend vers l'infini la proportion des suites $A^+$ est égale à la proportion des suites $S^+$ ceci se traduit par :

(6.1) $$\lim_{n \to +\infty} r_n(S^+) = \lim_{n \to +\infty} r_n(A^+)$$

Avec :

$r_n(S^+)$ La proportion des suites finies de type $S^+$

$r_n(A^+)$ La proportion des suites finies de type $A^+$

Equivalent à :
$$r_a(S^+) = r_a(A^+)$$
Ou encore
$$r_a(S^-) = r_a(A^-)$$

D'après le théorème 8.2 établi dans l'article [4] la proportion $r_a(A^+)$ est quasi-nulle

(6.2) $$r_a(A^+) = \lim_{n \to +\infty} r_n(A^+) = 0 \Rightarrow r_a(S^+) = 0$$

Sur la figure suivante en représente ces deux proportions en diagramme circulaire:

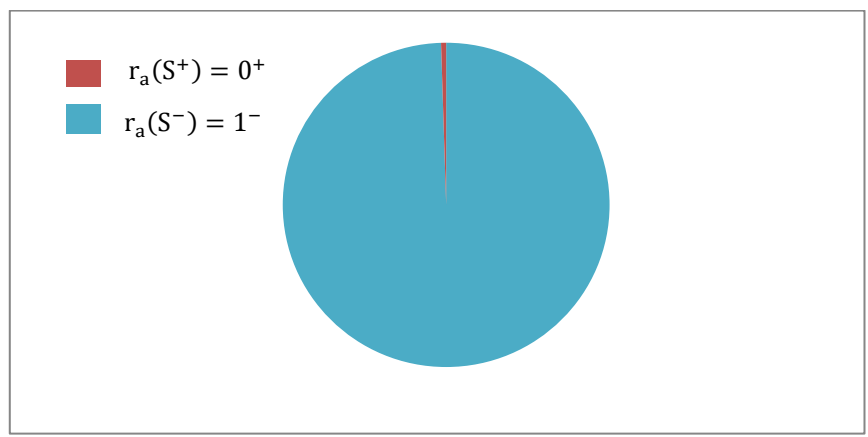

Figure 19 : Répartition des proportions des suites selon la monotonie

Donc la quasi-totalité des suites de Collatz de longueurs infinies sont décroissantes ou plus précisément admettant une limite à l' infini inférieure au son point de départ.

## 7.  Discussion des résultats obtenus

Cette étude nous permet de faire la distinction entre différentes catégories des suites de Collatz de longueurs infinies. Lorsque l'ordre dimensionnel d'un tableau complet de Collatz tend vers l'infini les différentes suites se comportement différemment. Une de répartition **possible** des suites infinies de Collatz dans le tableau complet d'ordre infini est représentée sur la figure suivante :

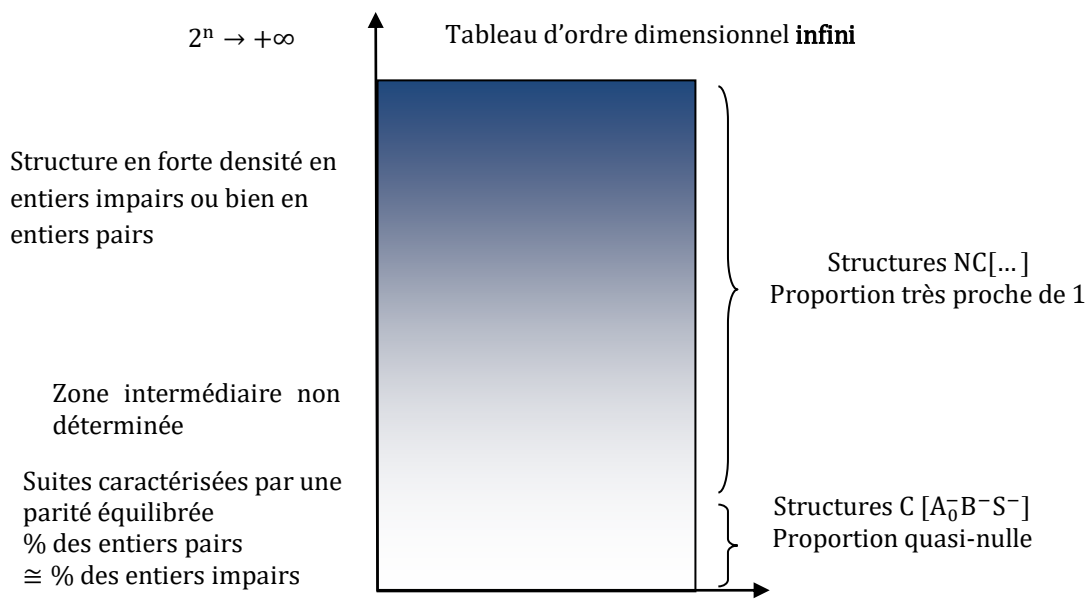

Figure 20: Répartition **apparente** de différentes catégories des suites infinies de Collatz dans un tableau complet d'ordre dimensionnel infini

Les suites qui vérifient la conjecture de Collatz sont toutes de la catégorie $R[A_0^- B^- S^-]$ alors que les autres suites n'admettent pas des termes finis et leurs vecteurs de parité correspondent à des suites binaires infinies **non convertibles** en suites de Collatz.

Cette conclusion ne peut être vraie que sous les conditions suivantes :

(1) Les ensembles $R[S^+ A_0^- B^-], R[S^+ A_0^- B^+], R[S^+ A_1^- B^+]$ et $R[S^+ A^+ B^+]$ doivent être tous des ensembles vides dans ces cas les seules suites infinies réalisables sont de type $R\, S^- A_0^- B^-$.

(2) La deuxième condition est le non existence d'un cycle stable autre le cycle (1,2).

Une fois ces deux conditions sont démontrées, on peut conclure que la conjecture de Collatz est vérifiée pour n'importe quelle suite de Collatz de premier terme un entier naturel non nul quelconque sans aucune exception.

Les travaux en cours portent sur la démonstration de l'impossibilité de l'existence des suites de Collatz appartenant aux catégories suivantes: $R\, S^+ A_0^- B^-, RS^+ A_0^- B^+, RS^+ A_1^- B^+$ et $RS^+ A^+ B^+$.

### *Loi de modération structurelle et équilibre de parité des suites infinies de Collatz*

La proportion des suites $A^+$ est quasi-nulle et la quasi-totalité des suites de Collatz de longueurs infinies sont de type $A^-$. Une suite de Collatz ne peut pas être dans cette catégorie ($[A^-]$) qui si la proportion des entiers impairs est strictement inférieure à 63%. De plus, la forte proportion en suites $A_0^-$ peut être expliquée par l'effet que la forte proportion des suites de Collatz contient une proportion en entiers impairs presque égale à celle en entiers pairs. En réalité cet équilibre de parité ou plutôt cette tendance vers un équilibre de parité peut être interprétée comme une loi de modération structurelle qui caractérise les suites des Collatz. Cette loi (si elle est vraie) signifie que les suites de Collatz tendent à égaliser leurs proportions en entiers pairs et impairs lorsque leurs longueurs tendent vers l'infini. Autrement :

Longueur de la suite $\to \infty \Rightarrow$ % des entiers pairs $\approx$ %des entiers impairs $\approx 50\%$

Cette dernière loi n'est pas encore démontrée mais la tendance de la proportion des suites $A^-$ vers 1 nous laisse croire à son existence et à sa vérité.

# Références